\documentclass[12pt]{amsart}
\usepackage{fullpage,amsfonts,a
msmath,amscd,amssymb,latexsym}
\input{xy}
\xyoption{all}
\xyoption{curve}
\title[subregular and simple]{Subregular representations of $\sl_n$
and simple singularities of type $A_{n-1}$}
\author{Iain Gordon}
\author{Dmitriy Rumynin}
\date{\today}
\newtheorem*{theorem}{Theorem}
\newtheorem*{prop}{Proposition}
\newtheorem*{lemma}{Lemma}

\newtheorem*{cor}{Corollary}

\newcommand{\CP}{\mathbb{C}}
\renewcommand{\sl}{{\mathfrak s}{\mathfrak l}}

\newcommand{\FF}{\mathcal{F}}

\renewcommand{\b}{{\mathfrak b}}

\newcommand{\bc}{B_{\chi,\lambda}}

\newcommand{\K}{{\mathbb K}}
\newcommand{\F}{{\mathbb F}}
\newcommand{\LF}{{\mathbb L}}
\newcommand{\ZK}{{\mathbb Z}/k{\mathbb Z}}
\newcommand{\T}{{\mathfrak T}}

\newcommand{\uc}{U_\chi}
\newcommand{\ucl}{U_{\chi,\lambda}}

\newcommand{\spr}{{\mathcal B}_\chi}
\newcommand{\ver}[1]{V(#1, \lambda)}

\newcommand{\gn}[2]{N_{#1}(#2)\gmd}

\newcommand{\uclt}{\ucl\text{-}T_0}
\newcommand{\UT}{\uc\text{-}T_0}
\newcommand{\cl}{C_{\LF}(n)}

\newcommand{\BT}{\bc\text{-}T_0}

\newcommand{\oo}{\overline}
\newcommand{\ra}{\rightarrow}

\newcommand{\Hom}{{\mbox{\text Hom}}}
\newcommand{\End}{{\mbox{\text End}}}
\newcommand{\Rep}{{\mbox{\text Rep}}}
\newcommand{\Coh}{{\mbox{\text Coh}}}

\DeclareMathOperator{\lie}{Lie}

\DeclareMathOperator{\soc}{soc}

\DeclareMathOperator{\md}{-mod}

\DeclareMathOperator{\gmd}{-grmod}

\DeclareMathOperator{\Gmd}{-Grmod}

\DeclareMathOperator{\rad}{rad}

\DeclareMathOperator{\ext}{Ext}

\DeclareMathOperator{\im}{im}

\DeclareMathOperator{\ed}{\End}

\begin{document}

\begin{abstract}
Alexander Premet has stated the following problem: what is a relation
between subregular nilpotent representations of a classical
semisimple restricted Lie algebra and non-commutative deformations
of the corresponding singularities? We solve this problem for type
$A$. Using the McKay correspondence, we relate the solution to
bases in equivariant $K$-theory introduced by Lusztig.
\end{abstract}

\thanks{Both authors were visiting and partially supported by MSRI while this research was begun and extend their thanks to that institution.
Much of the research for this paper was undertaken while the first
author was supported by TMR grant ERB FMRX-CT97-0100 at the
University of Bielefeld. The second author was partially supported
by EPSRC. }

\maketitle
\section{Introduction}
\subsection{}
Fix an integer $n>0$. Let $\K$ (respectively $\LF$) be an
algebraically closed field of characteristic $p$ (respectively
arbitrary characteristic). We will assume throughout that $p$ and
$n$ are coprime. Let $\F$ be the prime subfield of $\K$. Given an
integer $m\in{\mathbb Z}$ we denote its reduction modulo $p$ by
$\oo{m}\in\F$.

\subsection{}
Quantisations of Kleinian singularities were introduced by Hodges
\cite{hod} for type $A$ and later by Crawley-Boevey and Holland
\cite{cra} for all types. Premet has suggested a possible
relationship between these quantisations and subregular
representations of a simple Lie algebra. In \cite{pre6} Premet
examines this in characteristic zero.

\subsection{}
The first half of the paper is concerned with exploring Premet's
suggestion in the modular type $A$ case. We refer the reader to
later sections of the paper for unexplained definitions and
notation. Let $U$ be the enveloping algebra of
$\mathfrak{sl}_n(\K)$ and for $v\in \K[t]$ let $T(v)$ be Hodges'
non-commutative deformation of a Kleinian singularity of type $A$
over $\K$. Let $U_{\chi}$ be the reduced enveloping algebra for a
subregular nilpotent functional $\chi$ and let $U_{\chi,\lambda}$
be the central reduction of $U_{\chi}$ determined by the weight
$\lambda$. We show in Theorem \ref{prethm} that, for $p>n$, there
exists a polynomial $v_{\lambda}\in \K[t]$ such that
$U_{\chi,\lambda} \cong \mbox{Mat} (t(v_{\lambda}))$, where
$t(v_{\lambda})$ is a central reduction of $T(v_{\lambda})$.
Moreover this isomorphism respects a natural $\mathbb{Z}$-grading
on the algebras $U_{\chi,\lambda}$ and $t(v_{\lambda})$ and so
induces an equivalence between the category of $\uclt$-modules and
a category of suitably graded $t(v_{\lambda})$-modules.

\subsection{}
A crucial step in the above theorem is the comparison of
$U_{\chi,\lambda}$ and $t(v_{\lambda})$ with a basic algebra we
call the \textit{no-cycle algebra}. This algebra is defined over
any field $\LF$, depends on a positive integer $k$, and is denoted
$N_{\LF}(k)$. The no-cycle algebra is a string algebra and hence
its indecomposable modules can be described by a simple
combinatorial procedure. As an application of this comparison we
give a sufficient condition for the indecomposability of a baby
Verma module belonging to the block of $U_{\chi}$ determined by a
regular weight $\lambda$. This is formulated in terms of the
geometry of the Springer fibre $\mathcal{B}_{\chi}$.

\subsection{}
\label{intcat} The second half of the paper concerns conjectures
of Lusztig on a relationship between equivariant $K$-theory and
modular representations of simple Lie algebras, \cite{lusbas2}.
We work only with subregular representations of type A,
where these conjectures have been confirmed numerically by
Jantzen, \cite{ja1}. Let ${\bf H}$ be the minimal
desingularisation over $\CP$ of a Kleinian singularity of type $A$
and let ${\bf H}_0$ be the corresponding exceptional divisor. Both
${\bf H}$ and ${\bf H_0}$ admit an action of a two dimensional
torus, $T$. Lusztig has constructed an action of an affine
Iwahori-Hecke algebra on both $K_T({\bf H})$ and $K_T({\bf H}_0)$,
the $T$-equivariant $K$-groups of ${\bf H}$ and ${\bf H_0}$
respectively. Using Kapranov and Vasserot's interpretation of the
McKay correspondence in terms of derived categories, \cite{K-V},
we show in Theorem \ref{cat} that the above action on $K_T({\bf
H}_0)$ admits a categorification in $D_T({\bf H}_0)$, the bounded
derived category of $T$-equivariant sheaves on ${\bf H}_0$. This
essentially means that $D_T({\bf H}_0)$ admits an action of the
affine braid group which specialises to Lusztig's Iwahori-Hecke
algebra action on $K_T({\bf H}_0)$. We also describe a form and a
duality on $D_T({\bf H}_0)$ which specialise to those considered
by Lusztig on $K_T({\bf H}_0)$ in \cite{lusbas1}. The reader
should beware that we do not concern ourselves with the existence
of a braid group action as described in \cite{del}; in the
literature, our action is sometimes referred to as a weak action.

\subsection{}
The two themes of the paper are linked by the no-cycle algebra,
$N_{\LF}(n)$. There is a natural bigrading on $N_{\LF}(n)$. In the
characteristic zero case there is a triangulated functor from
$D(N_{\CP}(n)\gmd)$, the bounded derived category of bigraded
$N_{\CP}(n)$-modules, to $D_T({\bf H}_0)$. This functor induces an
isomorphism on $K$-theory, sending simple $N_{\CP}(n)$-modules to
the signed basis of $K_T({\bf H}_0)$ constructed by Lusztig in
\cite{lusnot}. In the modular case, forgetting half of the
bigrading on $N_{\K}(n)$, there is an isomorphism between the
$K$-groups of graded $N_{\K}(n)$-modules and of $\uclt$-modules
for regular $\lambda$. Under this isomorphism the operations of
wall-crossing functors on $U_{\chi}$-modules correspond to the
generators of an Iwahori-Hecke algebra action, a characteristic
$p$ analogue of that in \ref{intcat}. The category of bigraded
$N_{\K}(n)$-modules is a mixed category lying over the usual
category of $\uclt$-modules.

\subsection{}
There are several clear directions for future work arising from
this paper. Firstly, the hypothesis $p>n$ should be weakened to
$p$ and $n$ being coprime. It would also be highly desirable to
extend the results to other types. In Section \ref{finalr} we
indicate how to extend our results to type $B$. It seems likely,
however, that the techniques used here are not sufficient for this
in general. Furthermore, we expect $N_{\CP}(n)$-modules to
correspond to a central reduction of the subregular
representations of a quantum group of type $A$ at a root of unity.
We have not, however, checked this here.

\subsection{}
The paper is organised as follows. In Section 2 we introduce the
notation we will require when we deal with categories having group
actions. We define the no-cycle algebra in Section 3 and describe
some of its indecomposable representations. In Section 4 we study
subregular representations of $\mathfrak{sl}_n(\K)$ and relate
them to the no-cycle algebra. Section 5 takes a brief look at
Gr\"{o}bner-Shirshov bases for non-commutative algebras and their
representations. In Section 6 we study Hodges' quantisation using,
in particular, the results of the previous section. In Section 7
we prove  Theorem \ref{prethm}, whilst in Section 8 we present the
application to baby Verma modules. We remind the reader about
equivariant sheaves and Hilbert schemes in Section 9. In Section
10 we find a new interpretation of Lusztig's signed basis (in the
subregular situation) and in Section 11 we prove Theorem
\ref{cat}. Section 12 sees our final remarks, tying together the
two themes.

\section{$G$-categories}
\subsection{}
\label{defgr} Let $\mathcal{C}$ be an abelian $\LF$-category and
let $G$ be a group acting on $\mathcal{C}$. So, for every $g\in
G$, we have an exact \textit{shift functor}, $[g]$, together with
natural isomorphisms $[g]\circ [g'] \ra [gg']$. Sometimes in the
literature such actions are called ``weak'' as opposed to
``strong'' actions, which satisfy associativity constraints for
natural isomorphisms. We do not use the term ``weak'' since we are
not interested in associativity constraints.

We call $\mathcal{C}$  a $G$-category if an action of  $G$ on
$\mathcal{C}$ is fixed. If $\mathcal{C}$ and $\mathcal{D}$ are
both $G$-categories, we say that the functor $\Phi :\mathcal{C}\ra
\mathcal{D}$ is a $G$-functor if
functors $\Phi \circ [g]$ and $[g] \circ \Phi$
are naturally isomorphic for every $g\in G$.

A $G$-functor $\Phi:\mathcal{C}\ra \mathcal{D}$ is a
$G$-\textit{equivalence} if there exists a $G$-functor
$\Psi:\mathcal{D}\ra\mathcal{C}$ such that $\Psi\Phi \cong
1_{\mathcal{C}}$ and $\Phi\Psi \cong 1_{\mathcal{D}}$. We say
$\mathcal{C}$ and $\mathcal{D}$ are equivalent $G$-categories.
Note that an equivalence that is a $G$-graded functor need not be
a  $G$-equivalence \cite[Section 5]{gorgre}.

\subsection{}
\label{grrg} Let $R=\oplus_{g\in G}R_g$ be a $G$-graded
(noetherian) algebra, that is $R_gR_{g'} \subseteq R_{gg'}$. A
$G$-graded $R$-module is an $R$-module, $M$, together with a
$K$-space decomposition $M=\oplus_{g\in G}M_g$ satisfying
$R_g \cdot M_{g'} \subseteq M_{gg'}$. The category of $G$-graded
(finitely generated) $R$-modules, denoted $R\Gmd$ ($R\gmd$) is an
example of a $G$-category. By definition, we have
$
(M[g])_{g'} = M_{g'g^{-1}}.
$
for all $g,g'\in G$.

If $G$ acts on an algebra $R$ by automorphisms then $R\md$ is
another example of a $G$-category. The group $G$ acts by twisting.

\subsection{}
\label{grhom}
Given $X,Y \in \mathcal{C}$ and $g\in G$, 
set $\Hom(X,Y)_g = \Hom_{\mathcal{C}}(X[g], Y)$. We define
\[
\Hom(X, Y) = \bigoplus_{g\in G} \Hom(X,Y)_g,
\]
a $G$-graded vector space. The identification
$\Hom_{\mathcal{C}}(X[g], Y) \cong \Hom_{\mathcal{C}}(X[gg'],
Y[g'])$ yields a composition law for $X,Y,Z\in\mathcal{C}$:
$\Hom(Y,Z)_{g'} \times \Hom(X,Y)_{g} \ra \Hom(X,Z)_{gg'}$. Then,
for the examples appearing in \ref{grrg}, the space $\End(X) =
\Hom(X,X)$ becomes a $G$-graded $K$-algebra and $\Hom(X,Y)$ a
$G$-graded $\End(X)$-module. The functor $\Hom(X, -)$ is a
$G$-graded functor from $\mathcal{C}$ to $\End(X)\Gmd$.

\subsection{}
\label{grder}
%Let $D(\mathcal{C})$ denote the bounded derived category
%of $\mathcal{C}$.
{\em A triangulated $G$-category} is a triangulated category
with exact shift functors $[g]$,
 $g\in G$, and  natural isomorphisms
$[g]\circ [g'] \ra [gg']$.
The bounded derived category $D(\mathcal{C})$
of a $G$-category
is  a triangulated $G$-category.
Analogous to \ref{defgr}, we have the notions of
triangulated $G$-functors and equivalences of
triangulated $G$-categories.

\subsection{}
\label{grK}
Let $K(\mathcal{C})$ denote the Grothendieck group of $\mathcal{C}$. The Grothendieck group of $D(\mathcal{C})$, denoted $K'(\mathcal{C})$, is defined to be the free abelian group generated by the isomorphism classes of objects
of $D(\mathcal{C})$ subject to the relations
$[M]- [M']-[M'']$
for every distinguished triangle
$$
M' \longrightarrow M \longrightarrow M'' \longrightarrow M' [1] \leadsto .
$$
Taking the Euler characteristic induces an isomorphism 
${\bf i}: K'(\mathcal{C}) \rightarrow K(\mathcal{C})$.

Let $\mathcal C$ be a $G$-category.
Since shift functors are exact, they induce an action of $\mathbb{Z}[G]$, the group algebra of $G$, on both $K(\mathcal{C})$ and $K'(\mathcal{C})$, which commutes with ${\bf i}$. 
%the isomorphism between them. 
From now on, we will identify $K(\mathcal{C})$ and $K'(\mathcal{C})$.

\section{The no-cycle algebra}
\label{nc}
\subsection{}
\label{quiv} Let $k\in \mathbb{N}$ be a fixed integer. Let $Q$ be
the directed graph with $k$ vertices and $2k$ edges labelled $a_i$
and $b_i$ for $i\in \ZK$, see Figure 1.
\begin{figure}
\[
\begin{picture}(200,155)(10,20)
\put(82,153){$i$}
\put(86,144){\circle*{4}}
\put(45,148){$a_{i-1}$}
\put(58,118){$b_{i-1}$}
\put(115,145){$a_i$}
\put(100,118){$b_i$}
\put(175,118){$a_{i+1}$}
\put(140,95){$b_{i+1}$}
\put(94,146){\vector(3,-1){40}}
\put(82,140){\vector(-3,-1){40}}
\put(129,126){\vector(-3,1){40}}
\put(40,135){\vector(3,1){40}}
\put(36,126){\circle*{4}}
\put(5,124){$i-1$}
\put(138,126){\circle*{4}}
\put(141,133){$i+1$}
\put(146,126){\vector(3,-2){40}}
\put(181,91){\vector(-3,2){40}}
\put(188,92){\circle*{4}}
\put(195,91){$i+2$}
\qbezier[25](36,110)(40,0)(178,80)
\end{picture}
\]
\caption{}
\end{figure}
Let $\LF Q$  be the path algebra of $Q$.
The \textit{no-cycle algebra}, $N_{\LF}(k)$, is the quotient of $\LF Q$ by the two sided ideal generated by all non-trivial paths in $Q$ which start and end at the same vertex. By inspection, $N_{\LF}(k)$ is a $k(2k-1)$-dimensional algebra.
\subsection{}
\label{grnc} The algebra $N_{\LF}(k)$ is $\mathbb{Z}$-graded: a
vertex idempotent $e_i$ has degree $0$, $a_i$ has degree $-1$ and
$b_i$ has degree $1$. Following \ref{grrg} the category of
$\mathbb{Z}$-graded modules will be denoted by $\gn{\LF}{k}$.
\subsection{}
We need some notation before describing several
$N_{\LF}(k)$-modules. For $i\in\ZK$, introduce formal inverses of
the arrows $a_i$ (respectively $b_i$), written $a_i^{\ast}$
(respectively $b_i^{\ast}$). The head (respectively tail) of an
arrow, $c$, is denoted $h(c)$ (respectively $t(c)$), and we define
$h(c^{\ast}) = t(c)$ (respectively $t(c^{\ast}) = h(c)$). We form
\textit{formal paths} of length t, $c_1\ldots c_t$, where each
$c_j$ is of the form $c$ or $c^\ast$ for some arrow $c$ and
$t(c_j) = h(c_{j+1})$. Given a formal path $c_1\ldots c_t$, we
define its inverse to be $c_t^\ast\ldots c_1^\ast$, where
$(c^\ast)^\ast$ equals $c$.

\subsection{}
For $t\leq k$, let $S_t$ be the set of formal paths $c_1\ldots
c_t$ such that if $c_j = a_i$ (respectively $b_i$) then $c_{j+1}$
is either $a_{i-1}$ or $b_{i-1}^\ast$ (respectively $b_{i+1}$ or
$a_{i+1}^\ast$), and similarly if $c_j = a_i^\ast$ (respectively
$b_i^\ast$). Furthermore if $t=k$ then exclude from $S_k$ the
formal paths consisting entirely of $a$'s or entirely of $b$'s and
also the inverses of such formal paths. For $t < k$ (respectively
$t=k$) let $\rho_t$ be the equivalence relation on $S_t$ which
identifies a formal path with its inverse (respectively its cyclic
permutations and their inverses). Let $W_t$ be a set of
equivalence class representatives in $S_t$ for $\rho_t$.

\subsection{}
\label{string} For $t<k$, an element $C=c_1\ldots c_t\in W_t$
defines a $t+1$-dimensional indecomposable {\em string}
$N_{\LF}(k)$-module $St(C)$. A basis is given by $\{ z_0,\ldots
,z_t \}$ where, for $1\leq j\leq t$, the element $z_j$ is
concentrated at vertex $h(c_j)$, and $z_0$ is concentrated at
vertex $t(c_1)$. For $1\leq j \leq t$ if $c_j = c$, an arrow,
define $c(z_j) = z_{j-1}$, whilst if $c_j = c^\ast$, define
$c(z_{j-1}) = z_j$. All other arrows in $Q$ act as zero.

\subsection{}
\label{band} For $t=k$, and element $C= c_1\ldots c_k \in W_k$
together with a scalar $\lambda\in \LF^*$, define an
indecomposable {\em band} $N_{\LF}(k)$-module $Bd_\lambda (C)$. A
basis is given by $\{ z_0, \ldots z_{k-1} \}$ where, for $0\leq
j\leq k-1$, the element $z_j$ is concentrated at vertex $h(c_j)$.
If $c_k=c$, an arrow, define $c(z_0) = \lambda z_{k-1}$, whilst if
$c_k= c^\ast$, define $c(z_{k-1}) = \lambda^{-1} z_0$. For $1\leq
j \leq k-1$, if $c_j = c$, an arrow, then $c(z_j) = z_{j-1}$,
whilst if $c_j = c^\ast$ define $c(z_{j-1}) = z_j$. All other
arrows of $Q$ act as zero.

\subsection{}
The algebra $N_{\LF}(k)$ belongs to a family of tame algebras
called \textit{string algebras}. Any non-zero indecomposable
$N_{\LF}(k)$-module of dimension no greater than $k$ is isomorphic
to either $St(C)$ or $Bd_\lambda (C)$ for some unique $C$ and
$\lambda$, see \cite[Section 3]{butrin}. Observe that the modules
$St(C)$ admit  a $\mathbb{Z}$-grading compatible with the
$\mathbb{Z}$-grading on $N_{\LF}(k)$ introduced in \ref{grnc}.

\section{The reduced enveloping algebra}
\label{red}
\subsection{} Let $\chi\in \mathfrak{sl}_{n}(\LF)^\ast$ be
the functional vanishing on upper triangular matrices
and defined as follows on the strictly lower triangular matrices
\[
\chi (E_{i,j}) = \begin{cases}
1 \quad &\text{if $i=j+1$ and $1\leq j\leq n-2$}, \\
0 &\text{otherwise}.
\end{cases}
\]
Let $P$ (respectively $P^+$) be the weight lattice (respectively
dominant weights) of $SL_n(\LF)$ with respect to the standard
choice of torus and Borel subgroup. Let $\{ \varpi_1,\ldots
,\varpi_{n-1} \}$ be the fundamental weights and let $\rho =
\varpi_1 +\cdots + \varpi_{n-1}$. The Weyl group $W$ is the
symmetric group $\mathfrak{S}_{n}$ acting  on $P$. The $W$-orbit
through $\lambda\in P$ contains a unique representative in $P^+$.
We also need a dot action given by $$ w \bullet \lambda = w
(\lambda + \rho) - \rho,
\
w\in W, \lambda \in P.
$$

\subsection{}
\label{springerfibre} Let $\mathcal{B}$ be the flag variety. As a
set this consists of all Borel subalgebras of
$\mathfrak{sl}_{n}(\LF)$, that is those subalgebras which are
conjugate under $SL_{n}(\LF)$ to the upper triangular matrices.
The cotangent bundle of $\mathcal{B}$ is naturally identified with
the variety $$ \tilde{\mathcal{N}} = \{ (x, {\mathfrak b}) : x
({\mathfrak b}) = 0 \} \subset \mathfrak{sl}_{n}(\LF)^\ast \times
\mathcal{B} $$ where the first projection $\pi_1$ becomes the
moment map. The Springer fibre $\spr$ is the subvariety
$\pi_1^{-1} (\chi)$ of $\mathcal{B}$.

There is a simple way to parametrise $\spr$, \cite[Section
6.3]{slo}. Given a basis $u_i$ of $\LF^n$, let  $\FF(u_{1},
\ldots, u_n)$ be a flag with the span of $u_1,\  \ldots \ , u_k$
as the $k$-dimensional space. Let $v_i$ be an element of the
standard basis of $\LF^n$ so that $E_{i,j}v_j =v_i$. We introduce
the flag $$ \FF_{k,\alpha} = \FF (v_1, v_2,\; \ldots \; , v_{k-1},
v_{k} + \alpha v_n, v_n, v_{k+1}, v_{k+2}, \; \ldots \; ,v_{n-1})
$$ for all $(k,\alpha)\in (\{1,\ldots,n-1\}\times \LF)\cup (0,0)$.
The irreducible components of $\mathcal{B}_\chi$ are projective
lines $\Pi_k$, $1\leq k \leq n-1$ where $$ \Pi_k = \{
\FF_{n-k,\alpha} \; | \; \alpha \in \LF\} \cup \{ \FF_{n-k-1,
0}\}. $$ For $2\leq k \leq n-1$ the components $\Pi_{k-1}$ and
$\Pi_{k}$ intersect transversally at a point $p_{k-1,k}=\FF_{n-k,
0}$. Components $\Pi_i$ and $\Pi_j$ with $|i-j|>1$ do not
intersect.

Consider the following one parameter subgroup of the diagonal
matrices in $SL_{n}(\LF)$
\begin{equation*}
%\label{torus1}
T_0 = \{ \nu (\tau) = \tau E_{1,1} + \tau E_{2,2} \cdots + \tau
E_{n-1,n-1} + \tau^{1-n} E_{n,n}: \tau \in \LF^*\}.
\end{equation*}
Notice that $T_0$ stabilises $\mathcal{B}_{\chi}$, since $\nu
(\tau) \cdot \FF_{i,\alpha} = \FF_{i,\tau^{-n}\alpha}$.
\subsection{}
Let us further assume that $\LF = \CP$. By the Jacobson-Morozov
theorem there exists an $\mathfrak{sl}_2$-triple $e,h,f\in
\mathfrak{sl}_n(\CP)$ such that Tr$(ex)=\chi (x)$ for each $x\in
\mathfrak{sl}_n(\CP)$. Let $\mathcal{N}$ be the variety of
nilpotent elements in $\mathfrak{sl}_n(\CP)$. Let $$ V_\chi = \{
\mu \in \mathfrak{sl}_n(\CP)^\ast \; | \; \forall x \in
\mathfrak{sl}_n(\CP) \ \mu ([x,f]) = \chi ([x,f]) \}. $$ By
\cite[Theorem 6.4 and Section 7.4]{slo} $V_\chi$ is a Kleinian
singularity of type $A_{n-1}$ and
$\Lambda_{\chi}=\pi_1^{-1}(V_\chi)$ is its minimal
desingularisation with the exceptional fibre $\mathcal{B}_\chi$.

\subsection{}
{\it For the rest of this section we will assume that $p>n$}. We
expect, however, all results to hold under the weaker condition
that $p$ and $n$ are coprime. The proof of Proposition \ref{pro}
is the crucial point where we require the restriction $p>n$ and
all other results of this section, up to \ref{lienoc} inclusive,
will continue to hold if this proposition can be proved under the
weaker hypothesis.
\subsection{}
We are going to work with the Lie algebra $\mathfrak{sl}_{n}(\K)$
now. For any element $X\in \mathfrak{sl}_{n}(\K)$ let $X^{[p]}\in
\mathfrak{sl}_{n}(\K)$ denote the $p$-th power of $X$.

Let $U(\mathfrak{sl}_{n}(\K))$ be the universal enveloping algebra
of $\mathfrak{sl}_{n}(\K)$. For any $X\in \mathfrak{sl}_{n}(\K)$
the element $X^p - X^{[p]}\in U(\mathfrak{sl}_{n}(\K))$ is
central. We will study representations of the following
\textit{subregular reduced enveloping algebra}
\[
\uc := \frac{U(\mathfrak{sl}_{n}(\K))}{(X^p - X^{[p]} - \chi(X)^p : X\in \mathfrak{sl}_{n}(\K))}.
\]

Let $\Lambda = P/pP$, an $\F$-vector space, and let $\psi :P \ra
\Lambda$ be the quotient map. Both the usual and the dot actions
of $W$ on $P$ pass to $\Lambda$. Let $W(\lambda )$ be  the
stabiliser of $\lambda\in \Lambda$ under the dot action. Set
\[
C_0 = \{ \lambda\in P: \lambda+ \rho= \sum r_i\varpi_i \text{ with
}
 r_i \geq 0 \text{ and } p \geq (r_1+\cdots +r_{n-1}) \}.
\]
A weight in the interior of $C_0$ is called {\it regular}.

\subsection{Blocks}
\label{blocks}
By \cite[Theorem 3.18]{brogor1} we have a block decomposition
\[
\uc = \bigoplus_{\lambda} \bc,
\]
where $\lambda$ runs over a set of representatives of the
$W\bullet$-orbits on $\Lambda$.

Let ${\mathfrak h}$ be the diagonal Cartan subalgebra of
$\mathfrak{sl}_{n}(\K)$. The Weyl group acts by algebra
automorphisms on the polynomial ring $S({\mathfrak h})$. For
$\lambda\in\Lambda$ the partial coinvariants give a local, graded
algebra
\[
C_{\lambda} = \frac{ S({\mathfrak h})^{W(\lambda)}} {
(S({\mathfrak h})_+^W) }.
\]
Thanks to \cite[Theorem 10]{mirrum} and \cite[Theorem 8.2]{pre6} there is an injective algebra homomorphism $C_{\lambda} \longrightarrow \bc$, whose image is central in $\bc$. Henceforth, we will identify $C_{\lambda}$ with its image in $\bc$.
\subsection{}
We will be concerned with the category of finite dimensional
$\uc$-modules, $\uc\md$, or more specifically the subcategory
of $\bc$-modules. These categories have graded
analogues which we introduce now.

By construction $\chi (tXt^{-1}) = \chi (X)$ for all
$X\in\mathfrak{sl}_{n}(\K)$ and $t\in T_0$. As a result the action
of $T_0$ on $U(\mathfrak{sl}_{n}(\K))$ passes to an action on the
quotient $\uc$.

Following Jantzen \cite{ja3},
a $\UT$-module to be a finite dimensional vector space $V$ over
$K$ that has a structure both as a $\uc$-module and as a rational
$T_0$-module such that the following compatibility conditions hold:
\begin{enumerate}
\item
We have $t(Xv) = (tXt^{-1})tv$ for all $X\in \mathfrak{sl}_{n}(\K),t\in T_0$ and $v\in V$;
\item
The restriction of the $\mathfrak{sl}_{n}(\K)$-action on $V$ to $\lie(T_0)$ is equal to the derivative of the $T_0$-action on $V$.
\end{enumerate}
We obtain the category $\UT\md$, whose objects are the
$\UT$-modules and whose morphisms are the $T_0$-equivariant
$U_{\chi}$-module homomorphisms.

For $i\in\mathbb{Z}$, there are shift functors $[i]: \UT\md
\longrightarrow \UT\md$. These send a given $\UT$-module $V$ to the
object having the same $U_{\chi}$-module structure but with $T_0$
acting by $\nu(\tau).v = \tau^{ip} \nu(\tau)v$ for $\nu(\tau)\in T_0$ and all $v\in V$. This makes $\uclt\md$  a $\mathbb{Z}$-category.

By \cite[9.3]{gorpre} the full $\mathbb{Z}$-subcategory $\BT\md$
of $\UT\md$ is well-defined. Its objects are $\bc$-modules with a
compatible rational $T_0$-action. The projection functor
$\text{pr}_{\lambda} : \UT\md \longrightarrow \BT\md$ is a
$\mathbb{Z}$-functor.
\subsection{}
\label{grad}
Let $F: \UT\md\longrightarrow U_{\chi}\md$ denote the functor which
forgets the $T_0$-structure. The objects of $U_{\chi}\md$ which are in
the image of $F$ are called \textit{gradable}. It follows from \cite[Corollary 3.4]{gorgre} and
\cite[Corollary 1.4.1]{janmod} that the simple $U_{\chi}$-modules and
their projective covers are gradable. Moreover any lift of a simple
$U_{\chi}$-module is simple in $\UT\md$ and any lift of a projective
indecomposable $U_{\chi}$-module is projective indecomposable in $\UT\md$. Suppose $M$ is gradable,
that is there exists a $\UT$-module $V$ such that $F(V)=M$. Then,
by \cite[Remark 1.5]{janmod}, we have $F( \soc V) =\soc M$ and
$F(\rad V) = \rad M$.

For any $M,N\in\UT\md$, using the notation of \ref{grhom}, we have $\Hom_{U_{\chi}}(F(M), F(N))= \oplus_i \Hom_{\UT}(M[i],N)$, \cite[Section 2]{gorgre}.
\subsection{}
The category $\UT\md$ admits a contravariant self-equivalence, $D$
whose square is canonically isomorphic to the
identity functor, \cite[Sections 1.13 and 1.14]{jan}. Moreover $D$ fixes the simple modules in $\BT\md$, \cite[Proposition 2.16]{jan}.

\subsection{Translation functors}
Using the map $\psi: P\ra \Lambda$,  we will abuse notation by
writing $\BT\md$ and $\text{pr}_{\lambda}$ for $\lambda\in P$ (we
should really take the image of $\lambda$ under $\psi$). Given
$\lambda, \mu \in C_0$ we define a translation functor
\[
T_{\lambda}^{\mu}: \BT\md \ra B_{\chi,\mu}\text{-}T_0\md
\]
by $T_{\lambda}^{\mu}(V) = pr_{\mu}(E\otimes V)$, where $E$ is the
simple $SL_{n}(\K)$-module with the highest weight $w (\mu -
\lambda)\in P^+$ for some $w\in W$.

Note that we get (in general) more than one functor
$\bc\text{-}T_0\md \ra B_{\chi,\mu}\text{-}T_0\md$ for fixed $\lambda$ and $\mu$: if $\mu$ and $\mu'$ are two distinct weights in $C_0$ with $\pi(\mu)$ and $\pi(\mu')$ in the same
$W$-orbit then $T_{\lambda}^{\mu}$ and $T_{\lambda}^{\mu'}$ will be two (in general) distinct functors from $\bc\text{-}T_0\md$ to $B_{\chi,\mu}\text{-}T_0\md = B_{\chi,\mu'}\text{-}T_0\md$.

\subsection{Baby Verma modules}
 Given $\b \in \spr$ and $\lambda\in P$ we have a one dimensional representation of $U_0 (\b)$, the subalgebra of $\uc$ generated by the elements of $\b$, described as follows. Let $g\in SL_{n}(\K)$ be such that $g$ conjugates $\b$
to the upper triangular matrices in $\mathfrak{sl}_{n}(\K)$, say $\b_+$.
There is a one dimensional $U_0(\b_+)$-module which is annihilated by strictly
upper triangular matrices and on which diagonal
matrices act via $\psi(\lambda)$.
Conjugation by $g$ provides an isomorphism between $U_0(\b)$ and $U_0(\b_+)$, and so gives a one dimensional $U_0(\b)$-module, say $\K_{\lambda}$. It can be checked that this module is independent of the choice of $g\in SL_{n}(\K)$. Induction yields a \textit{baby Verma module}
\[
\ver{\b} = \uc \otimes_{U_0 (b)} \K_{\lambda}
\]
This is a $\bc$-module on which $C_{\lambda}$ acts by scalar multiplication.

\subsection{}
\label{uni} The module $\ver{\b_+}$ can be given the structure of
a $\BT$-module where $T_0$ acts on $1\otimes 1$ through $\lambda$.
Set $\ver{\b_+}' = D(\ver{\b_+})$ (it follows from
\cite[11.16(1)]{ja3} that this is a baby Verma module with respect
to the Borel subalgebra obtained by conjugating $\b_+$ by the
Coxeter element $s_1s_2 \ldots s_{n-1}$).

\begin{prop}{\cite[Theorem 2.6]{jan}}
Suppose $\lambda \in C_0$ with $\lambda +\rho = \sum r_i\varpi_i$
and let $r_0 = p- (r_1+ \cdots +r_{n-1})$.

(i) The category $\BT\md$ has simple modules $L_0, L_1, \ldots , L_{n-1}$ (up to isomorphism and shift) where the module $L_i$ has dimension $p^{(n^2-n-2)/2}r_{n-1-i}$ (if $r_{n-1-i}=0$ then $L_i$ should be omitted from the list of simple modules).

(ii) For each $i$ between $0$ and $n-1$ there exists a unserial
module $V_i\in \BT\md$ whose Loewy layers are $L_i, L_{i+1}[-1],
\ldots , L_{i-1}[1-n]$ (omit $V_i$ and $L_i$ whenever $r_{n-1-i} =
0$).

(iii)  For each $i$ between $0$ and $n-1$ there exists a unserial
module $V'_i\in \BT\md$ whose Loewy layers are $L_i, L_{i-1}[1],
\ldots ,L_{i+1}[n-1]$ (omit $V'_i$ and $L_i$ whenever $r_{n-1-i}=
0$).

(iv) For any $\mu\in X$ such that $\pi(\mu)$ and $\pi(\lambda)$
are in the same $W\bullet$-orbit there exists a unique $i$ and
$j\in \mathbb{Z}$ (respectively $i',j'$) such that $V(\b_+, \mu)
\cong V_i[j]$ (respectively $V(\b_+, \mu)' \cong V'_{i'}[j']$).
\end{prop}
\subsection{Endomorphisms}
\label{pro} Let $\lambda\in P$ and let $Z_{\lambda}$ be the centre
of $B_{\chi,\lambda}$. Recall $C_{\lambda}\subseteq Z_{\lambda}$.
Let $M$ be in $\BT\md$. There is a homomorphism
\[
\theta_M : Z_{\lambda} \longrightarrow \ed_{\bc}(F(M)),
\]
sending $z\in Z_{\lambda}$ to the endomorphism $(m \longmapsto
z\cdot m)$.
\begin{prop}
Let $\lambda \in C_0$ and $Q$ be a projective indecomposable
module in $\BT\md$. The homomorphism
\[
\theta_{Q}:Z_{\lambda} \longrightarrow \mbox{\em End}_{\bc}(F(Q))
\]
is surjective.
\end{prop}
\begin{proof}
To make the paper as self-contained as possible, we will give a
direct proof of this for regular weights $\lambda$, since the
general case relies on an unpublished result of Jantzen,
\cite{ja2}.

It is enough to prove this in the ungraded case. Let $P_0, \ldots
P_r$ be all distinct (up to isomorphism) projective indecomposable
$B_{\chi,\lambda}$-modules. The algebra $C_{\lambda}$ has a simple
socle \cite[Corollary 3.9]{gor}. Therefore if
Ann$_{C_{\lambda}}(P_i)$ is non-zero it must contain the socle of
$C_{\lambda}$. The equality
$0=\mbox{Ann}_{C_{\lambda}}(B_{\chi,\lambda})
=\cap_i\;\mbox{Ann}_{C_{\lambda}}(P_i)$, implies that
$C_{\lambda}$ acts faithfully on at least one projective
indecomposable, say $P_0$. By Proposition \ref{uni} and a result
of Jantzen \cite[Proposition 10.11]{ja3} the dimension of
$\ed_{B_{\chi,\lambda}}(P_0)=\dim C_{\lambda}$, so $Z_{\lambda}$
generates $\ed_{B_{\chi,\lambda}}(P_0)$.

Now assume $\lambda$ is regular. It follows from \cite[Section
2.3]{jan} that we can find a translation functor $T$ such that
$T(P_0)\cong P_i$ and, by \cite[Section 11.21]{ja3}, that $T$ is a
self-equivalence of $B_{\chi,\lambda}\md$. Since $Z_{\lambda}$ can
be identified with the endomorphism ring of the identity functor
on $B_{\chi,\lambda}\md$, conjugation by $T$ induces a ring
automorphism of $Z_{\lambda}$, say $\tilde{T}$. It follows that
$\ed_{B_{\chi,\lambda}}(P_i)$ is generated by
$\tilde{T}(C_{\lambda}) \subseteq Z_{\lambda}$, as claimed.

If $\lambda$ is not regular the above argument fails since it need
no longer be true that we can find a translation functor $T$ which
is a self-equivalence and sends $P_0$ to $P_i$. In this situation
we use the following fact, \cite[C.6 Claim 2]{ja2}: if the highest
weight of $P_i$ (in the graded category) belongs to $C_0$ and does
not lie on the affine wall, then the canonical map $C_{\lambda}
\ra \ed_{B_{\chi,\lambda}}(P_i)$ is an isomorphism. Hence it is
enough to show that we can find a representative of the
isomorphism class of $P_i$ whose highest weight belongs to $C_0$
and does not lie on the affine wall. A straightforward calculation
shows that this follows from \cite[Section 2.3]{ja1}.
\end{proof}
\subsection{}
It is possible to {\it slightly} weaken the hypothesis $p>n$ when
using the results of \cite{ja2}, and hence the hypothesis of this
whole section. The proof of Jantzen's results, however, do not
apparently generalise to the best case where $p$ and $n$ are
coprime. We leave these details to the interested reader.
\subsection{}
\label{nocf} Let $J$ be the unique maximal ideal of $Z_{\lambda}$.
\begin{lemma}
Let $P_i$ be the projective cover of $L_i$ in $\BT\md$. Then $[P_i
/JP_i
:L_i[j]]=\delta_{j0}$.
\end{lemma}
\begin{proof}
It is enough to prove this for ungraded $\bc$-modules. Suppose
that $F(L_i)$ is a composition factor of $F(P_i)/JF(P_i)$. Thus
$F(L_i)$ appears as a direct summand of $\rad^m F(P_i) /
\rad^{m+1} F(P_i)$ for some $m\in \mathbb{N}$. Hence we have a
commutative diagram
\[
\xymatrix{ F(P_i) \ar@{->}[r] \ar@{-->}[drr] &
\rad^mF(P_i)/\rad^{m+1}F(P_i) \ar@{->}[r] &
F(P_i)/\rad^{m+1}F(P_i) \\ & & F(P_i) \ar@{->}[u] }
\]
Thanks to Lemma \ref{pro} there exists $z\in Z_{\lambda}$ such
that the above endomorphism of $F(P_i)$ is multiplication by $z$.
Then $z\notin J$ since, by hypothesis, the composition factor
$F(L_i)$ does not lie in $JF(P_i)$. Since $Z_{\lambda}$ is local
it follows that the endomorphism is an isomorphism and so $F(L_i)$
lies in the head of $F(P_i)$ as required.
\end{proof}
\subsection{Two central reductions}
Let $\tilde{J}$ be the unique maximal ideal of $C_{\lambda}$ and
recall that $J$ is the unique maximal ideal of $Z_{\lambda}$. We
introduce two central reductions
\[
\tilde{U}_{\chi,\lambda}  = \frac{\bc}{\tilde{J}\bc}, \ \ \ \ucl =
\frac{B_{\chi,\lambda}}{JB_{\chi,\lambda}}.
\]
Since $\tilde{J}$ and $J$ are homogeneous both
$\tilde{U}_{\chi,\lambda}$ and $\ucl$ inherit
$\mathbb{Z}$-gradings from $B_{\chi,\lambda}$. The category of
graded modules $\uclt\md$ is thus a full subcategory of $\BT\md$
and there is $\mathbb{Z}$-action on $\uclt\md$, inherited from
$\BT\md$.
\subsection{} \label{lienoc}
The main result of this section follows.
\begin{prop}
Suppose $\lambda\in C_0$ with $\lambda+\rho = \sum r_i \varpi_i$,
and let $r_0 = p-(r_1+\ldots +r_{n-1})$. Let $k$ be the number of
non-zero $r_i$'s. Then $\ucl$ is Morita equivalent to the no-cycle
algebra $N_{\K}(k)$. Moreover, if $\lambda$ is regular there is a
$\mathbb{Z}$-equivalence of categories
\[
\uclt\md \longrightarrow \gn{\K}{n}.
\]
\end{prop}
\begin{proof}
Let $A$ be a finite dimensional algebra with simple modules $S_1, \ldots , S_r$. The Gabriel quiver of $A$ is the directed graph with vertices labelled from $1$ to $r$ and $\dim \ext_{A}^1(S_i,S_j)$ edges from $i$ to $j$. By \cite[Proposition 4.17]{ben} $A$ is Morita equivalent to the the path algebra of its Gabriel quiver factored by some admissible ideal, that is an ideal generated by linear combinations of paths of length at least two.

Let $0\leq i_1
< i_2 < \ldots < i_k \leq n-1$ be such that $L_{i_t}\neq 0$, or,
equivalently, $r_{n-1-i_t}\neq 0$. Set $s_t = i_{t+1} - i_t$ for
$1\leq t < k$ and set $s_k = n + i_1 - i_k$. Since $L_i$ appears
only once as a composition factor of $V_i$, $Z(B_{\chi,\lambda})$
acts by scalars on $V_i$, making $V_i$ a $\ucl$-module. By
\cite[Proposition 2.19]{jan}, for $t\neq t'$ modulo $k$ and $j\in
\mathbb{Z}$
\[
\ext_{\uclt}^1 (L_{i_t}[j], L_{i_{t'}}) =
\left\{ \begin{array}{cl}
\K & \mbox{ if } t' = t +1, j=s_t \mbox{ or } t'=t-1, j= -s_{t-1}, \\
0 &  \mbox{ otherwise}
\end{array}
\right.
\]
Thus the Gabriel quiver of $\ucl$ is of the form \ref{quiv}, possibly with loops added at the vertices. Let $B$ be the quotient of this quiver which is Morita equivalent to $\ucl$. We will show $B$ is isomorphic to $N_{\K}(k)$.

The projective covers of the simple $B$-modules are spanned by the paths ending in a fixed vertex. Hence, Lemma \ref{nocf} shows that there can be no loops at vertices and further, that $B$ is therefore a quotient of $N_{\K}(k)$. In particular, its dimension is at most $k(2k-1)$.

Let $T_{i_t}$ be the kernel of the sum of two projections $V_{i_t}
\oplus V_{i_t}^\prime \rightarrow L_{i_t}$. By Proposition
\ref{uni}
\[
[F(T_{i_t}) : F(L_{i_{t'}})] =
\begin{cases}
2 \quad & \mbox{ if } \ t \neq t' \\
1 &  \mbox{ if } \ t = t'.
\end{cases}
\]
Since $F(T_{i_t})$ is a quotient module of $F(P_{i_t})$,
the dimension of $B$ can be estimated by
\[
\dim B = \ed (\bigoplus_{t=1}^k F(P_{i_t}))
= \sum_{t,t'} [F(P_{i_t}) :F(L_{i_{t'}})] \geq  \sum_{t,t'} [F(T_{i_t}) :F(L_{i_{t'}})] =
k(2k-1).
\]
We deduce that $B\cong N_{\K}(k)$, proving the first statement of the theorem. This also proves that $T_{i_t}$ is the projective cover of $L_{i_t}$ in $\uclt\md$.

Let $\lambda$ be regular, so that $k=n$, and let $T = \oplus T_{i_t}$. Thanks to \ref{grhom} and \ref{grad} the algebra
$E=\ed_{\ucl}(F(T))^{\mbox{\tiny op}}$ has a $\mathbb{Z}$-grading.
By \cite[Theorem 5.4]{gorgre} $\uclt\md$ and $N_{\K}(n)\md$ are equivalent $\mathbb{Z}$-categories if $E \cong N_{\K}(n)$ as a $\mathbb{Z}$-graded algebra. But, up to a choice of scalars, $b_i$ corresponds to a $\uclt$-module homomorphism sending $T_i[1]$ to $T_{i+1}$, and $a_i$ corresponds to the a $\uclt$-module homomorphism sending $T_{i+1}[-1]$ to $T_i$. This proves the second claim.
\end{proof}

\subsection{Two central reductions (II)}
\label{rem} In general the inclusion $C_{\lambda}\subseteq
Z_{\lambda}$ is strict, \cite[Corollary 3.9]{gor}. However, the
natural homomorphism $\tilde{U}_{\chi,\lambda}\ra\ucl$ is an
isomorphism. This follows from the fact in the proof of
Proposition \ref{pro} that the map $$ \theta_1: C_{\lambda}
\longrightarrow \ed_{\bc}(F(P)) $$ is an isomorphism for any
projective indecomposable module $P$. The arguments of \ref{nocf}
and \ref{lienoc} are then valid for $\tilde{U}_{\chi,\lambda}$,
from which the isomorphism follows. We expect this continues to
hold under the weaker hypothesis that $p$ and $n$ are coprime.

\section{Gr\"{o}bner-Shirshov Bases}
\subsection{}We are going to use Gr\"{o}bner-Shirshov bases for
associative algebras (see \cite{bok} for two-sided ideals and
\cite{kan} for one-sided ideals). In this section we quickly
explain the technique to make our paper self-contained. Although
the version for one-sided ideals \cite{kan} is sufficient for our
ends, we generalise to arbitrary modules to avoid repetitions.

\subsection{}
Let $R=\LF <X_1,X_2,\ldots,X_l>$ be a free associative algebra.
Let $F$ be a free $R$-module with generators $Y_1, \ldots , Y_k$.
Although $l$ and $k$ are natural numbers here, one can use, with
certain care,  the technique for transfinite ordinals.

Let $\mathcal R$ be the set of monomials in $R$,
$\mathcal F$ the set of monomials in $F$.
The set ${\mathcal R} \cup {\mathcal F}$
admits a partial multiplication with a two-sided unit
$1_R$ (one agrees that $m1=m$ for $m \in {\mathcal F}$).
We always make the assumption
that a product is defined when we write the product.
We start with a linear order $\succ$ on
${\mathcal R} \cup {\mathcal F}$
such that
\begin{itemize}
\item $\forall z \in {\mathcal R} \cup {\mathcal F} \ \ \ z \succ 1_R$;
\item $z_1 \succ z_2 \ \Longrightarrow
wz_1 v\succ w z_2 v$;
\item $\forall z\in {\mathcal R}\cup {\mathcal F} $ the set
$\{ w \in {\mathcal R}\cup {\mathcal F}| z \succ w \}$ is finite.
\end{itemize}
A degree lexicographical order is most practical
but there are different orders.
%A classification
%of all such orders is an interesting problem.
For a non-zero element $f$ of $R\cup F$
we denote the highest term of $f$
by $\oo{f}$.

\subsection{}
A pair of subsets $S\subseteq R$ and $T\subseteq F$ determine an
algebra $A=R/RSR$ and a left $A$-module $M= A\otimes_R (F/RT)=
F/(RT + RSF)$. The technique of Gr\"{o}bner-Shirshov pairs allows
to solve questions about $M$ by producing an explicit basis of
$M$.

Each $f\in S\cup T$ gives
a rewriting rule $\oo{f}\rightarrow \oo{f}-f$.
We write $a\leadsto b$ if
$b$ can be obtained from $a$ by using
rewriting rules.
Note that one cannot rewrite a monomial $\oo{f}t$
unless $t=1$ or $f\in S$.

\subsection{Composition}
For certain $f,g \in R\cup F$, $w\in {\mathcal R} \cup {\mathcal F}$,
we can form a composition $(f,g)_w$.
If $w=\oo{f}V  = W \oo{g}= WZV$ for some
$W,Z,V \in{\mathcal R} \cup {\mathcal F}$
with $Z \neq 1$ then the composition is
$$
(f,g)_w = fV - W g .
$$
If $w=W \oo{f}V  = \oo{g}$ for some
$W,V \in{\mathcal R} \cup {\mathcal F}$
then the composition is
$$
(f,g)_w = W fV - g.
$$
These two cases are mutually exclusive.

A pair $(S,T)$ is a Gr\"{o}bner-Shirshov pair if $(f,g)_w\leadsto
0$ for all possible $f,g \in S\cup T$ and $w \in{\mathcal R} \cup
{\mathcal F}$.

The following version of Shirshov's composition lemma can be
proved by standard methods \cite{bok}. If $k=1$ and $T=\emptyset$
then the statement is the standard version of Shirshov's
composition lemma \cite{bok}. If $k=1$ and $T$ arbitrary then it
is a version for left ideals \cite{kan}.

\subsection{Shirshov's composition lemma}
{\em For $S$ and $T$ as above, we consider the set of monomials $$
B = \{ Z \in {\mathcal F}\;  | \; \forall f\in S\cup T \  \forall
W, V \ Z \neq W \oo{f} V \}. $$ If $(S,T)$ is a
Gr\"{o}bner-Shirshov pair then the image of $B$ is  a basis of $M$
as an $\LF$-vector space.

Moreover, for every $(S,T)$ there exits a Gr\"{o}bner-Shirshov
pair $(S^\prime ,T^\prime)$ such that $RSR = RS^\prime R$ and $RSF
+RT =RS^\prime F +R T^\prime$. }

\subsection{Buchberger's algorithm}
A proof of existence of a
Gr\"{o}bner-Shirshov pair uses transfinite recursion,
called Buchberger's algorithm. 
It proceeds
as follows. One starts with $(S_0,T_0)=(S,T)$. Given $(S_m,T_m)$ we
produce the next pair $(S_{m+1},T_{m+1})$ such that $RS_mR =
RS_{m+1} R$ and $RS_mF +RT_m =RS_{m+1} F +R T_{m+1}$. Consider all
possible compostions $(f,g)_w$ with $f,g\in S_m\cup T_m$. To each
such composition, apply a sequence of rewriting rules
$\oo{v}\rightarrow \oo{v}-v$ with $v\in S_m\cup T_m$ so that
$(f,g)_w \leadsto [f,g]_w$ and $[f,g]_w$ cannot be rewritten any
further.
Note that the element $[f,g]_w$ is not canonical since
we choose a sequence of rewriting rules to use.
Another sequence can give a different answer.
Define the following sets
$$ I_m =\{[f,g]_w \ | \ f,g \in S_m\cup T_m, \ [f,g]_w \neq 0\}, $$
$$ J_m =\{g \ | \ f,g \in
S_m\cup T_m,\  (f,g)_w = WfV-g,\ [f,g]_w \neq 0\}. $$ Now we can
make the recursion step, $$ S_{m+1} = (S_m \cup (I_m \cap
{\mathcal R}))\setminus J_m, \ T_{m+1} = (T_m \cup (I_m \cap
{\mathcal F}))\setminus J_m. $$

\subsection{Termination of Buchberger's algorithm}
If $(S,T)$ is a Gr\"{o}bner-Shirshov pair then $S_1=S_0$,
$T_1=T_0$, and the procedure terminates immediately.
%In many applications
%the procedure terminates after several steps.
If $S$ is a Gr\"{o}bner-Shirshov basis (equivalently
$(S,\emptyset)$ is a Gr\"{o}bner-Shirshov pair), and $R/RSR$ is
noetherian then the procedure terminates after finitely many
steps.

\section{Hodges' Quantisation}
\label{hq}
\subsection{Kleinian singularities} \label{ks} Let
$\zeta\in \K$ be a primitive root of unity of degree $n$. Set
\[
\Gamma = \left\{ g^i : g = \left( \begin{matrix} \zeta  & 0 \\ 0 &
\zeta^{-1}
\end{matrix}
\right)
\right\},
\]
a subgroup of $SL_2(\K)$. The natural action of $\Gamma$ on $\K^2$
induces an action on $\K [X, Y]$: $g\cdot X = \zeta X$, $g\cdot Y
= \zeta^{-1}Y$. The invariants of $\K [X,Y]$ under this action are
generated by $X^{n}, XY$ and $Y^{n}$. Thus, the orbit space
$\K^2/\Gamma$ has co-ordinate ring
\[
\mathcal{O}(\K^2/\Gamma) = \K [X^{n}, XY, Y^{n}] \cong \frac{\K [
A,B,H]}{(AB- H^{n})}.
\]
The variety $\K^2/\Gamma$ has an isolated singularity at $0$, a
\textit{Kleinian singularity of type $A_{n-1}$}.

\subsection{}
Let $v(z)\in \K [z]$ be a polynomial of degree $n$, whose roots lie in $\F$. Following \cite{hod}, we define an associative algebra, $T(v)$, over $\K$ with generators $a,b$ and $h$ satisfying the relations
\begin{equation}
\label{TV}
ha = a(h+1), hb = b(h-1), ba = v(h), ab=v(h-1).
\end{equation}
There exists a filtration on $T(v)$
such that
$
\mbox{gr}(T(v))\cong \K[A,B,H]/
(AB-H^{n}).
$
In other words,
$T(v)$ is a deformation of a Kleinian singularity of type $A_{n-1}$.

Using the translation $h\mapsto h-1$ we will assume without loss of generality that $0$ is a root of $v(z)$.

\subsection{A $\mathbb{Z}$-category}
There is also a $\mathbb{Z}$-grading on $T(v)$: we assign $a$
degree $n$, $b$ degree $-n$ and $h$ degree $0$. We want to study a
category of graded $T(v)$-modules similar in spirit to the
construction of $\BT\md$ in Section \ref{grad}.

Since $p$ and $n$ are coprime,
 we can find an inverse of  $\oo{n}$ in $\mathbb{F}$, say $\oo{c}$.
We will consider the full subcategory of finitely generated $\mathbb{Z}$-graded $T(v)$-modules consisting of objects
\[
M = \bigoplus_{j\in\mathbb{Z}} M_j,
\]
such that $h$ acts on $M_j$ through scalar multiplication by $\oo{j}\oo{c}$. 
Note that $a\cdot M_j \subseteq M_{j+n}$ 
(respectively $b\cdot M_j \subseteq M_{j-n}$) showing that 
this definition is compatible with the relation $ha=a(h+1)$ 
(respectively $hb= b(h-1)$). Denote this category by $T(v)\gmd$.

For $i\in\mathbb{Z}$, there is a shift functor $[i]: T(v)\gmd \ra T(v)\gmd$: given $M\in T(v)\gmd$, set $(M[i])_j = M_{j-pi}$ for all $j\in \mathbb{Z}$. This makes $T(v)\gmd$ a $\mathbb{Z}$-category.

\subsection{}
We will be interested in a finite dimensional central quotient of $T(v)$.
\begin{lemma}
The centre of $T(v)$ is generated by $a^p, b^p$ and $h^p-h$. It is isomorphic to the algebra of functions on a type $A_{n-1}$ Kleinian singularity.
\end{lemma}
\begin{proof}
It is straightforward to check that $a^p$, $b^p$ and $h^p - h$ are central elements. By construction $T(v)$ is a free $\K [h]$-module with basis $\{ a^i, b^j :i,j\geq 0\}$.
It follows from the relations in $T(v)$ that the degrees of the homogeneous components of any non-zero central elements must be a multiple of $pn$.
Since $a^p$ and $b^p$ are central we must find which polynomials in $h$ are central. Let $q(h)$ be such a polynomial. Since $aq(h) = q(h+1)a$ we deduce that the roots of $q$ are invariant under integer addition. It follows that $q(h)$ is a polynomial in $h^p-h$ as required.

Using the defining relations once more we have
\[
a^pb^p = v(h)v(h+1)\ldots v(h+p-1).
\]
Since $v(h)$ has degree $n$ it follows that $a^pb^p = (h^p - h)^{n}$. Hence, the centre of $T(v)$ is a quotient of the ring of functions of a Kleinian singularity of type $A_{n-1}$. Any proper quotient of the ring of functions on a Kleinian singularity has dimension 0 or 1. Thus since $T(v)$ is finitely generated over its centre and has Gelfand-Kirillov dimension 2, the centre must be the entire ring of functions.
\end{proof}

\subsection{}Now we can introduce the protagonist of this section:
\begin{equation*}
%\label{tv}
t(v)= \frac{T(v)}{(a^p, b^p, h^p-h)}.
\end{equation*}
Since the ideal $(a^p, b^p, h^p -h)$ is homogeneous, $t(v)$ inherits a $\mathbb{Z}$-grading from $T(v)$. We denote the full subcategory of $T(v)\gmd$ consisting of $\mathbb{Z}$-graded $t(v)$-modules by $t(v)\gmd$.

\subsection{}
In order to study $t(v)$ we introduce an intermediate algebra
\begin{equation}
\label{Tv}
\T(v)=T(v)/(a^p,b^p).
\end{equation}
Let us use a degree lexicographical order on non-commutative associative
monomials in $a$ of degree 1, $b$ of degree $2n-1$, and $h$
of degree 1 with $h > b > a$.

The relations of (\ref{TV}) already form a Gr\"{o}bner-Shirshov
basis.
It follows that monomials not containing $ha$, $hb$, $ba$, or $ab$
as submonomials form a basis of $T(v)$.
%In particular, this proves
%the claim
%of \cite{hod} that
%\begin{equation} \label{sum}
%T(v) = \ldots \oplus a^2 \K[h] \oplus a\K[h] \oplus \K[h]
%\oplus b \K[h] \ldots
%\end{equation}
%
%\marginpar{Do we need a restr. on $p$?}

For any polynomial $f(z)\in \K [z]$ and a positive integer $i$, we denote
$$
f_{(i)}(z) = \prod_{k=0}^{i-1} f(z+\oo{k}), \
f_{(-i)}(z) = \prod_{k=0}^{i-1} f(z-\oo{k}).
$$

\begin{lemma} \label{shirshov}
The following relations, together with those in (\ref{TV}) and
(\ref{Tv}), form a Gr\"{o}bner-Shirshov basis of $\T(v)$,
\begin{eqnarray}
a^{p-i} v_{(-i)}(h-1) & \mbox{ for } & i=1,\; \ldots \; ,p,
\\
b^{p-i} v_{(i)}(h)  & \mbox{ for } & i=1,\; \ldots \; ,p-1.
\end{eqnarray}
\end{lemma}
\begin{proof}
Let us obtain all relations in (\ref{shirshov})
recursively. For $i=1, \;\ldots\; p$,
$$(ba-v(h),a^{p-i}v_{(-i)}(h-1))_{ba^{p-i}h^{ni}}=
ba^{p-i}v_{(-i)}(h-1) - (ba-v(h))a^{p-i-1}h^{ni}=
$$
$$
= ba^{p-i}(v_{(-i)}(h-1)-h^{ni}) + v(h)a^{p-i-1}h^{ni}
\leadsto
v(h) a^{p-i-1}(v_{(-i)}(h-1)-h^{ni}) +
$$
$$
+ v(h)a^{p-i-1}h^{ni} =
v(h) a^{p-i-1}v_{(-i)}(h-1) \leadsto
a v(h+1) a^{p-i-2}v_{(-i)}(h-1) \leadsto
\ldots
$$
$$
\leadsto
a^{p-i-1} v(h+p-i-1) v_{(-i)}(h-1) =
a^{p-i-1}v_{(-i-1)}(h-1).
$$
Similarly, for $i= 1, \ldots , p-1$,
$
(b^{p-i} v_{(i)}(h), ab-v(h-1))_{ab^{p-i}h^{ni}}
\leadsto b^{p-i-1} v_{(i+1)}(h).
$
Now we need to show that all remaining compositions are trivial.
The highest terms of defining relations
are $ha$, $hb$, $ba$, $ab$, $a^{p-i}h^{ni}$, and $b^{p-i}h^{ni}$.
Let us make certain that all compositions are trivial.
Firstly,
$$            %composition 1
(ba-v(h),ab-v(h-1))_{bab}=
bv(h-1)-v(h)b
\leadsto bv(h-1)-bv(h-1)= 0.
$$
Similarly,
$
(ab-v(h-1),ba-v(h))_{aba}\leadsto 0.
$
Then
$$               %composition 2
(hb-bh+b,b^{p-i}v_{(i)}(h))_{hb^{p-i}h^{ni}} =
hb^{p-i}(h^{ni}- v_{(i)}(h)) -bhb^{p-i-1}h^{ni}
+b^{p-i}h^{ni}
\leadsto
$$
$$
\leadsto
b^{p-i}( (h-p+i)(h^{ni}- v_{(i)}(h)) -(h-p+i+1)h^{ni}
+h^{ni} =
b^{p-i} v_{(i)}(h) h
\leadsto 0.
$$
Analagously, the compositions
$
(a^{p-i}v_{(-i)}(h),ab-v(h-1))_{a^pb} $ and $
(b^{p-i}v_{(i)}(h),ab-v(h-1))_{ab^p} $ are trivial. Then $$
(ha-ah-a,ab-v(h-1))_{hab}= hv(h-1) -ahb-ab \leadsto hv(h-1) -
a(bh-b)-ab = $$ $$ =hv(h-1) -abh \leadsto 0. $$ Another possible
composition to consider is $$ (ha-ah-a,
a^{p-i}v_{(-i)}(h))_{ha^{p-i}h^{ni}} = ha^{p-i}(v_{(-i)}(h) -
h^{ni}) + aha^{p-i-1}h^{ni} +a^{p-i}h^{ni} \leadsto $$ $$ \leadsto
a^{p-i}( (h-i)(v_{(-i)}(h) - h^{ni}) + (h-i-1)h^{ni} +h^{ni}) =
a^{p-i}v_{(-i)}(h) (h-i) \leadsto 0. $$ The remaining compositions
$
(ha-ah-a, a^{p-i}h^{ni})_{a^{p-i}h^{ni}a}
$,
$
(hb-bh+b, a^{p-i}h^{ni})_{a^{p-i}h^{ni}b}$,
$
(ha-ah-a, b^{p-i}h^{ni})_{b^{p-i}h^{ni}a}$, and
$(hb, b^{p-i}h^{ni})_{b^{p-i}h^{ni}b}$
are trivial by a similar argument.
\end{proof}

\subsection{Corollary}
{\em The dimension of $\T(v)$ is $np^2$. Moreover,
there is a direct sum decomposition
\begin{equation} \label{sum2}
\T(v) =
[\bigoplus_{i=1}^{p} a^{p-i} \K[h]/(v_{(i)}(h-1))]
\oplus
[\bigoplus_{j=1}^{p-1} b^{j} \K[h]/(v_{(p-j)}(h))]
\end{equation}
}
\begin{proof}
The  direct sum decomposition~(\ref{sum2}) follows at once from
the description of the Gr\"{o}bner-Shirshov basis of $\T(v)$.
Adding dimensions of summands, we arrive at the dimension of
$\T(v)$, that is $2(n+2n+\ldots (p-1)n ) + pn = np^2 $.
\end{proof}

\subsection{}
\label{weedim}
If we write $\mbox{gcd}(f(Z),g(Z))$ for the greatest common
divisor of two polynomials then decomposition (\ref{sum2})
is inherited:
\begin{equation} \label{sum3}
t(v) =
\bigoplus_{i=1}^{p} \frac{a^{p-i} \K[h]}{(\mbox{\small gcd}
(v_{(i)}(h-1),h^p-h))}
\oplus
\bigoplus_{j=1}^{p-1} \frac{b^{j} \K[h]}{(\mbox{\small gcd}
(v_{(p-j)}(h),h^p-h))}.
\end{equation}
This decomposition allows one to compute the dimension of
$t(v)$.

Let $r_1,\ldots ,r_{n-1}\in\mathbb{Z}$ be such that $r_i\geq 0$ and $p\geq r_1+\ldots r_{n-1}$ and
\[
v(z) = \prod_{i=0}^{n-1} (z- (\oo{r_1}+\ldots +\oo{r_i})).
\]
Set $r_0 = p- (r_1+\ldots +r_{n-1})$.
\begin{cor}
The dimension of $t(v)$ is $2p^2-\sum_{i=0}^{n-1} r_i^2$.
\end{cor}
\begin{proof}
If $f(z,i)=z-(\oo{r_1}+\ldots +\oo{r_i})$ then the roots of $f_{(j)}(z,i)$ are $\oo{r_1}+\ldots \oo{r_i}, \ldots ,\oo{r_1}+\ldots +\oo{r_i} -\oo{j}$. If $j\geq r_i$ then $\oo{r_1}\ldots \oo{r_{i-1}}, \ldots ,\oo{r_1} +\ldots +\oo{r_i} - \oo{j}$ are already roots of $f_{(j)}(z,i+1)$. Thus, the dimension of the first summand in (\ref{sum3})
is
\[
p^2 - \sum_{i=0}^{n-1} (1+2+\ldots + (r_i-1)+r_i +r_i + \ldots
r_i),
\]
where the number of summands in the parenthesis is $p$. This sum equals
$p^2+(p-\sum_{i=0}^{n-1} r_i^2)/2$. Similarly, the second summand has dimension $p^2 - (p + \sum_{i=0}^{n-1} r_i^2)/2$, so that the total dimension
is $2p^2-\sum_{i=0}^{n-1} r_i^2$.
\end{proof}
\subsection{Baby Verma modules (II)}
\label{byV}
Let $u$ (respectively $u'$) be the subalgebra of $t(v)$ generated by $a$ and $h$ (respectively $b$ and $h$). Similarly let $U$ (respectively $U'$) be the subalgebra of $\T(v)$ generated by $a$ and $h$ (respectively $b$ and $h$).

We introduce two sets of baby Verma modules. For $\lambda\in \mathbb{Z}$ let $\K_{\lambda}$ (respectively $\K_{\lambda}'$) be the one dimensional $U$-module (respectively $U'$-module) with basis $|0>$ (respectively $|1>$) and
$$
h|0>=\oo{\lambda} |0>, \ \ a|0>=0,
\ \ \ \ \ \
h|1>=\oo{\lambda} |1>, \ \ b|1>=0.
$$
The baby Verma module $V(\lambda)$ (respectively $V(\lambda)'$) is
\[
V(\lambda) = \T(v) \otimes_U \K_{\lambda} \quad (\mbox{respectively }V(\lambda)' = \T(v)\otimes_{U'} \K_{\lambda}').
\]

\subsection{Lemma}{ \em
(i) If $\oo{\lambda}$ is not a root of $v(z)$ then
$V(\lambda) =0$, whilst if $\oo{\lambda}$ is a root of $v(z)$ then
$V(\lambda)$ has a basis of $p$ elements
$|0>, b|0>, \ldots, b^{p-1}|0>$. We have
\begin{equation*}
ab^k|0>=v(\oo{\lambda} -\oo{k})b^{k-1}|0>,  \ \
hb^k|0>=(\oo{\lambda} -\oo{k})b^k|0>.
\end{equation*}
(ii) If $\oo{\lambda} -1$ is not a root of $v(z)$ then
$V(\lambda)' =0$, whilst if $\oo{\lambda} -1$ is a root of $v(z)$ then
$V(\lambda)'$ has a basis of $p$ elements
$|1>, a|1>, \ldots, a^{p-1}|1>$. We have
\begin{equation*}
ba^k|1>=v(\oo{\lambda} +\oo{k}-1)a^{k-1}|1>,  \ \
ha^k|1>=(\oo{\lambda} +\oo{k})a^k|1>.
\end{equation*}
}
\begin{proof}
(i) For the generator
$|0>\in V(\lambda)$ one observes that
\[
0=ba|0>-(ba-v(h))|0>=v(h)|0>=v(\oo{\lambda})|0>.
\]
Thus if $\oo{\lambda}$ is not a root of $v(z)$ then
$|0>=0$ and $V(\lambda) =0$.

Let $\oo{\lambda}$ be a root of $v(z)$. Let $S$ be the
Gr\"{o}bner-Shirshov basis of $\T(v)$ constructed in
Lemma~\ref{shirshov}. The module $V(\lambda)$ is determined by the
pair $(S,\{(h-\oo{\lambda})|0>,\ a|0>\} )$ and this turns out to
be a Gr\"{o}bner-Shirshov pair. Indeed, there are three elements
in $S$ whose leading monomials end with $a$:
\[
(a|0>,ha-ah-a)_a = h a|0> - (ha-ah-a)|0>= (ah+a)|0> \leadsto
(\lambda a +a)|0> \leadsto 0,
\]
\[
(a|0>,ba-v (h))_a =b a|0> - (ba-v (h))|0>= v(h)|0> \leadsto
v(\oo{\lambda})|0> =0,
\]
\[
(a^p, a|0>)_a =a^p|0> - a^{p-1} a|0> =0.
\]
Elements of $S$
whose leading monomials end with $h$
fall into two types:
$$
(a^{p-i} v_{(i)}(h-1)|0>, (h- \oo{\lambda} )|0>)_h =
(\oo{\lambda} a^{p-i} h^{ni-1} +
a^{p-i} (v_{(i)}(h-1)-h^{ni}))|0> \leadsto
$$
$$
v_{(i)}(\oo{\lambda} -1) a^{p-i}|0> \leadsto 0,
$$
$$
(b^{p-i} v_{(i)}(h),(h- \oo{\lambda})|0>)_h =
(\oo{\lambda} b^{p-i} h^{ni-1} +
b^{p-i} (v_{(i)}(h)-h^{ni}))|0> \leadsto
b^{p-i} v_{(i)}(\oo{\lambda})|0> = 0.
$$
Direct computation now yields the formulas for the action.

(ii) The proof is analogous.
\end{proof}
Thanks to the lemma we can consider $V(\lambda)$ and $V(\lambda)'$ as objects in $T(v)\gmd$. Indeed if $V(\lambda)$ (respectively $V(\lambda)'$) is non-zero we let $b^k|0>$ (respectively $a^k|1>$) span the $(\lambda-k) n$ (respectively $(\lambda + k)n$) homogeneous component.

\subsection{}
\label{hoduni}
Under our assumptions, we have $(\oo{\lambda} -\oo{k})^p=\oo{\lambda} -\oo{k}$.
It follows that $(h^p-h)V(\lambda) =0$ (respectively $(h^p-h)V({\lambda})'=0$), and so $V(\lambda)$ and $V({\lambda})'$ are objects in $t(v)\gmd$.
\begin{prop}
Let $v(z) = \prod_{i=0}^{n-1} (z- (\oo{r_1}+\ldots +\oo{r_i}))$
be as in \ref{weedim} and let $r_0 = p -(r_1+\ldots +r_{n-1})$.

(i) The category $t(v)\gmd$ has simple modules $L_0, \ldots , L_{n-1}$ (up to isomorphism and shift) where the dimension of $L_i$ is $r_{n-1-i}$ (if $r_{n-1-i}=0$ then $L_i$ should be omitted from the list of simple modules).

(ii) For each $i$ lying between $0$ and $n-1$ there exists a uniserial module $V_i\in t(v)\gmd$ whose Loewy layers are $L_i, L_{i+1}[-1],\ldots, L_{i-1}[1-n]$ (omit $V_i$ and $L_i$ whenever $r_{n-1-i}=0$).

(iii)  For each $i$ lying between $0$ and $n-1$ there exists a unserial module $V'_i\in t(v)\gmd$ whose Loewy layers are $L_i, L_{i-1}[1], \ldots ,L_{i+1}[n-1]$ (omit $V_i'$ and $L_i$ whenever $r_{n-1-i}= 0$).

(iv) For any $\lambda\in \mathbb{Z}$ such that $\oo{\lambda}$ is a root of $v(z)$ there exists a unique $i$ and $j\in \mathbb{Z}$ (respectively $i',j'$) such that $V({\lambda}) \cong V_i[j]$ (respectively $V({\lambda + 1})' \cong V'_{i'}[j']$).
\end{prop}
\begin{proof}
Set $V_0 = V(0)$ and $V_i = V(r_1+\ldots +r_{n-1-i})[i]$. Note that if $r_{n-1-i}=0$ then $V_i\cong V_{i+1}[-1]$. Since every $t(v)$-module has a $u$-fixed point we see that any simple $t(v)$-module is a quotient of $F(V_i)$ for some $i$. By Lemma \ref{byV}(i) $F(V_i)$ is isomorphic to $\K [b]/(b^p)$ as a $\K [b]/(b^p)$-module. Since $\K [b]/(b^p)$ is a local algebra it follows that $V_i$ has a simple head.

In $V_i$ the element $b^{r_{n-1-i}}|0>$ is annihilated by $a$ and belongs to the $n(r_1+\ldots +r_{n-1-(i+1)}) + pi$ component, so there is a graded $t(v)$ homomorphism
\[
\theta_i : V_{i+1}[-1] \longrightarrow V_i.
\]
The cokernel of $\theta_i$, say $L_i$, has dimension $r_{n-1-i}$ and is simple if $r_{n-1-i}\neq 0$ since it is generated by any basis vector $b^i|0>$ it contains. If $r_{n-1-i}\neq 0$ it follows from Lemma \ref{hoduni} that $L_i$ has a unique $u$-fixed point, namely $|0>$. Hence if $r_{n-1-i}$ and $r_{n-1-j}$ are non-zero $L_i$ and $L_j$ are isomorphic if and only if $i=j$. This proves (i).

Since $V_i$ is simple-headed, $\dim V_i = p$ and $\sum_{i=0}^{n-1} r_i = p$ the chain of homomorphisms
\[
\begin{CD}
V_{i-1}[1-n] @> \theta_{i-2} >> V_{i-2}[2-n] @> \theta_{i-3} >> \cdots @>\theta_{i+1} >> V_{i+1}[-1] @> \theta_{i} >> V_i
\end{CD}
\]
proves (ii). The proof of (iii) is similar. Part (iv) is clear.
\end{proof}

\subsection{Proposition}
\label{hodnoc}
{\em
Let $v(z) = \prod_{i=0}^{n-1} (z- (\oo{r_1}+\ldots +\oo{r_i}))$ be as in \ref{weedim} and set $r_0 = p -(r_1+\ldots +r_{n-1})$.
Let $k$ be the number of non-zero $r_i$'s.
Then $t(v)$ is Morita equivalent to $N_{\K}(k)$.
Moreover, if $k=n$, there is a $\mathbb{Z}$-equivalence of categories
\[
t(v)\gmd \longrightarrow \gn{\K}{n}.
\] }
\begin{proof}
Let $0\leq i_1 \leq \ldots \leq i_k \leq n-1$ be such that $r_{n-1-i_t}\neq 0$. Let $Q_{i_t}$ be the projective cover of $L_{i_t}$ in $t(v)\gmd$. Recall the general formula, \cite[Section 1.7]{ben}
\[
\dim t(v) = \sum_{t=1}^k \dim Q_{i_t} \dim L_{i_t}.
\]
Let $T_{i_t}$ be the kernel of the sum of two projections $V_{i_t} \oplus V'_{i_t}\ra L_{i_t}$. Then $T_{i_t}$ has head isomorphic to $L_{i_t}$ so is a quotient of $Q_{i_t}$. By Lemma \ref{byV} and Proposition \ref{hoduni} $\dim T_{i_t} = 2p - r_{n-1-i_t}$. Using Lemma \ref{weedim} we find
\[
\sum_{t=1}^k \dim Q_{i_t} \dim L_{i_t} \geq \sum_{t=1}^k \dim
T_{i_t} \dim L_{i_t} =\sum_{t=1}^k (2p - r_{n-1-i_t})r_{n-1-i_t}
\]
\[
= 2p^2 - \sum_{t=1}^k r_{n-1-i_t}^2 = \dim t(v),
\]
proving that $T_{i_t} \cong Q_{i_t}$.

Let $T=\oplus T_{i_t}$. The basic algebra of $t(v)$ is
$\ed_{t(v)}(F(T))^{\mbox{\small{op}}}$. Let $b_t$ (respectively
$a_t$) be the homomorphism $F(T_{i_t}) \ra F(T_{i_{t+1}})$
(respectively $F(T_{i_{t+1}})\ra F(T_{i_t})$) associated to the
composition factor $F(L_{i_t})$ of $F(T_{i_{t+1}})$ (respectively
$F(T_{i_t})$) lying in the second Loewy layer of $F(V_{i_{t+1}})$
(respectively $F(V_{i_t}')$). It is straightforward to check that
$a_t$ and $b_t$, together with the idempotents arising from the
projections in $\mbox{End}_{t(v)}(F(T))$, generate the basic
algebra and satisfy the relations of the no-cycle algebra. Since
\[
\dim \ed_{t(v)}(F(T)) = k(2k-1)
\]
the first statement of the proposition follows.
The second statement is proved in the same manner as Proposition \ref{lienoc}.
\end{proof}

\section{Proof of Premet's conjecture}
\label{prec}
\subsection{}
We require $p>n$ for the following theorem. Thanks to \ref{rem} we
can replace $\ucl$ with $\tilde{U}_{\chi,\lambda}$ throughout, if
we wish.
\subsection{Theorem}
\label{prethm} { \em Let $p>n$. Suppose $\lambda \in C_0$ with
$\lambda+\rho = \sum r_i\varpi_i$ and let $v(z) =
\prod_{i=0}^{n-1} (z- (\oo{r_1}+\ldots +\oo{r_i}))$. There is an
isomorphism
\[
\ucl \cong \mbox{\em Mat}_{p^{(n^2-n-2)/2}}(t(v)).
\]
Moreover, there is a $\mathbb{Z}$-equivalence of categories
\[
\uclt\md \longrightarrow t(v)\gmd.
\] }
\begin{proof}
Set $r_0 = p -(r_1 +\ldots +r_{n-1})$ and let $k$ be the number of
non-zero $r_i$'s. Let $0\leq i_1\leq \ldots \leq i_k \leq n-1$ be
such that $r_{n-1-i_t}\neq 0$. Let $L_1,\ldots ,L_k$ (respectively
$M_1,\ldots , M_k$) be the simple $\uclt$-modules (respectively
graded $t(v)$-modules) appearing in Proposition \ref{uni}
(respectively Proposition \ref{hoduni}) and let $P_1,\ldots ,P_k$
(respectively $Q_1,\ldots ,Q_k$) be their projective covers. We
have
\[
\ucl \cong \ed_{\ucl}(\oplus F(P_t)^{p^{(n^2-n-2)/2}r_{n-1-i_t}})
\cong \mbox{Mat}_{p^{(n^2-n-2)/2}}(\ed_{\ucl}(\oplus F(P_t)^{r_{n-1-i_t}}))
\]
and
\[
t(v) \cong \ed_{t(v)}(\oplus F(Q_t)^{r_{n-1-i_t}}).
\]
Thanks to our construction of $P_t$ in Proposition \ref{lienoc} and $Q_t$ in Proposition \ref{hodnoc} we have a graded isomorphism
\[
\ed_{\ucl}(\oplus F(P_t)^{r_{n-1-i_t}}) \cong \ed_{t(v)}(\oplus F(Q_t)^{r_{n-1-i_t}}),
\]
proving the first statement of the theorem, together with an
equivalence. The equivalence is a ${\mathbb Z}$-equivalence
by \cite[Theorem 5.4]{gorgre}.
\end{proof}

\section{More on Baby Verma Modules}

\subsection{}
\label{deflem} We want to study baby Verma modules in ``general
position". To do so, we need a general lemma.
\begin{lemma}
Let $A$ be a finite dimensional $\LF$-algebra and $Y$ a connected
algebraic variety over $\LF$. Suppose $M_\alpha$, $\alpha\in Y$,
is a flat family of finite dimensional $A$-modules over $Y$. Then
the Grothendieck group element $[M_\alpha ]\in K_0(A)$ is
independent of $\alpha$.
\end{lemma}
\begin{proof} Let $B$ be the basic algebra of $A$. There exists a
$(B,A)$-bimodule $N$, flat over $A$, such that the functor
$N\otimes_A-$ induces an equivalence between the categories of
finite dimensional $A$-modules and finite dimensional $B$-modules.
Let $K_\alpha= N\otimes_A M_\alpha$, a flat family of $B$-modules
over $Y$. Given a primitive idempotent $e\in B$, it suffices to
show that the dimension of $eK_\alpha$ is independent of
$\alpha\in U$, a Zariski open neighbourhood of a point. Without
loss of generality we can trivialise the family $K\times U$. Then
$e$ defines an algebraic family of projection operators on the
finite dimensional vector space $K$, $$ e_\alpha: k \mapsto
\mbox{pr}_K (e\cdot (k,\alpha)).$$ Since the dimension of
$eK_\alpha$ is equal to the rank of $e_\alpha$, and the latter is
obviously constant, the lemma is proved.
\end{proof}
\begin{proof}(Alternative)
Without loss of generality we can assume that $Y$ is irreducible.
Let $n=\dim M_{\alpha}$. For any $a\in A$ let $c_a(\alpha)$ be the
characteristic polynomial corresponding to the action of $a$ on
$M_{\alpha}$. As this a polynomial of degree $n$ we obtain a
regular map $\chi_a : Y \ra K^n$ sending $\alpha$ to the
coefficients of $c_a(\alpha)$ (coefficients written in order $1,
t, t^2, \ldots$).

Given a simple $A$-module, $S$, let $d_s(\alpha) = \dim S
[M_{\alpha}:S]$. Thanks to \cite[Section 4, Theorem]{craboe} we
have
\[
d_s(\alpha) = \text{min}_{a\in Ann(S)} \{ \mbox{ord}_{t=0}
c_a(\alpha)\}.
\]
For any $0\leq m \leq n$ let $V_m \subseteq K^n$ be the closed set
defined by the vanishing of the first $m$ co-ordinates. Then
$\chi_a^{-1}(V_m)$ is closed and consists of the values $\alpha$
for which $\mbox{ord}_{t=0} c_a (\alpha) \geq m$. Hence, for $a\in
A$, there exists an open set, $U_a\subseteq Y$ such that for all
$\alpha\in U_a$ the value $\mbox{ord}_{t=0}c_a(\alpha)$ is
minimal. Let $m_a$ be this minimum.

The set $\{ m_a : a\in Ann(S) \}$ has a minimum, denoted $\mu$,
which is achieved at $a_{\mu}\in A$ say. Then for all $\alpha\in
U_{a_{\mu}}$ we have $d_S(\alpha) = a_{\mu}$. We deduce that the
set $\{ \alpha\in Y : d_S(\alpha ) \text{ minimal}\}$ is dense in
$K$ since it contains $U_{a_{\mu}}$.

Repeating the above for all simples we see there is a dense set in
$K$ where all the functions $d_S(\alpha)$ take minimal values. But
$\sum_S d_S(\alpha) = \dim M_{\alpha} = n$ is constant and
therefore they must take minimal values everywhere, proving the
lemma.
\end{proof}

\subsection{}
Given $\lambda \in X$, it is an interesting problem to describe the
isomorphism classes of all baby Verma modules $\ver{\b}$ as $\b$
runs over $\mathcal{B}_{\chi}$. If $\lambda$ is regular then the
description of $\spr$, Proposition \ref{uni} and Lemma
\ref{deflem} show that, for every $\b\in \spr$,
\[
[\ver{\b}] = \sum_{i=0}^{n-1} [L_i] \in K_0(U_{\chi,\lambda}).
\]
Therefore, on passing to the no-cycle algebra, we see that such
$\ver{\b}$ corresponds to a $N_{\K}(n)$-module $M(\b)$, such that
$e_iM(\b)$ is one dimensional for all $i$. Thanks to Section
\ref{nc}, all modules of this dimension are known.

Let $(k,\alpha) \in (\{1,\ldots, n-1\}\times \K\})\cup (0,0)$ and
let $\b_{k,\alpha}$ be the stabiliser of $\FF_{k,\alpha}$. Suppose
first that $\alpha=0$. Since the torus $T_0$ stabilises $\b_{k,0}$
for any $k$, twists by elements of $T_0$ provide a grading of
$\ver{\b_{k,0}}$. Therefore $M(\b_{k,\alpha})$ is gradable and, by
Section \ref{nc}, must be a direct sum of string modules if $n$ is
odd (note that some band modules are gradable for even $n$). When
$\lambda$ is regular, we expect that gradable band modules are not
baby Verma modules, and that baby Verma modules are
indecomposable.

The generic case is dealt with in the following proposition.

\subsection{{\bf Proposition}}
{\em 
Keep the above notation and let $\lambda\in X$ be a regular
weight, $1\leq k\leq n-1$ and $\alpha\neq 0$. Then the baby Verma
module $\ver{\b_{k,\alpha}}$ is indecomposable.
}
\begin{proof}
It follows from Section \ref{nc} that if $M(\b_{k,\alpha})$ is not
gradable it is necessarily a band module and hence indecomposable.
Thus, by Proposition \ref{lienoc}, it suffices to show that
$\ver{\b_{k,\alpha}}$ does not admit a $T_0$-grading.

Suppose for a contradiction that $\ver{\b_{k,\alpha}}$ admits a
$T_0$-grading. Let $L_k$ be the $\sl_2(\K)$-subalgebra generated
by $E_{n-k,n}, E_{n-k,n-k}-E_{n,n}$ and $E_{n,n-k}$. Any Borel
subalgebra belonging to $\Pi_k$ is uniquely determined by its
intersection with $L_k$. Let $\lambda_k =
\lambda(E_{n-k,n-k}-E_{n,n})$. A straightforward calculation shows
that the restriction of $\ver{\b_{k,\alpha}}$ to $L_k$ has a
direct summand isomorphic to the baby Verma module for $L_k$
induced from $\b_{k,\alpha}\cap L_k$ with highest weight
$\lambda_k$.

For $t\in T_0$ the element $t(1\otimes 1)\in \ver{\b_{k,\alpha}}$
is a highest weight vector for the Borel subalgebra
$t\cdot\b_{k,\alpha}$, yielding an isomorphism
$\ver{\b_{k,\alpha}} \cong \ver{t\cdot\b_{k,\alpha}}$. Since
$\lambda$ is regular $\lambda_k\neq -1$, and so, by \cite[Main
Theorem]{pre}, the baby Verma modules for $L_k$ induced from
different Borel subalgebras of $L_k$ with highest weight
$\lambda_k$ are not isomorphic. Hence, by the last paragraph, the
isomorphism $\ver{\b_{k,\alpha}}\cong \ver{t\cdot\b_{k,\alpha}}$
forces $\ver{\b_{k,\alpha}}$ to have infinitely many
non-isomorphic direct summands, a contradiction.
\end{proof}

\section{Equivariant $K$-theory and Hilbert schemes}
\subsection{Equivariant sheaves}
\label{def}
Let $G$ be an affine algebraic group acting rationally on a quasi-projective noetherian variety $X$. Let $\Coh_G(X)$ denote the abelian category of $G$-equivariant coherent sheaves on $X$, \cite[Chapter 5]{CG}. Let $D_G(X)$ denote the \textit{bounded} derived category of $\Coh_G(X)$.

\label{map}
Let $f:X \ra Y$ be a $G$-equivariant map. If $f$ is proper then the right derived functor of $f_{\ast}$ gives a pushforward
\[
{\bf R}f_{\ast} : D_G(X) \longrightarrow D_G(Y).
\]
Similarly, if $f$ has finite Tor dimension (for example if $f$ is flat or if $Y$ is smooth), then the left derived functor of $f^{\ast}$ gives a pullback
\[
{\bf L}f^{\ast} : D_G(Y) \longrightarrow D_G(X).
\]
\subsection{}
\label{sup}
Suppose $Z$ is a closed $G$-stable subvariety of $X$. Let $D_G(X, Z)$ denote the full subcategory of $D_G(X)$ whose objects are complexes with homology supported on $Z$. If $i : Z\ra X$ is inclusion then ${\bf R}i_{\ast} : D_G(Z) \ra D_G(X,Z)$ is an equivalence of categories.
\subsection{}
Let $K_G(X)$ denote the Grothendieck group of $\Coh_G(X)$, or equally, by \ref{grK}, of $D_G(X)$. The constructions of \ref{map} and \ref{sup} yield pushforward (respectively pullback) homomorphisms between the Grothendieck groups of $X$ and $Y$ for a proper (respectively finite Tor dimension) $G$-equivariant map $f:X\ra Y$ and, for any closed $G$-stable subvariety $Z$, an identification between $K_G(Z)$ and the Grothendieck group of complexes of equivariant sheaves on $X$ whose homology is supported on $Z$.
\label{easybun}

Let $f: X \ra \{pt \}$ be projection to a point. Since $\Coh_G(pt)$ can be identified with the category of finite dimensional $G$-modules, any $G$-module, say $M$, can be pulled back to a $G$-equivariant locally free sheaf on $X$, $f^{\ast}M$. Since the operation of tensoring by a locally free sheaf is exact we obtain a functor $D_G(X) \ra D_G(X)$ associated to each $G$-module and therefore an action of $K_G(pt) = \Rep(G)$, the representation ring of $G$, on $K_G(X)$.
\subsection{Equivariant line bundles on $\mathbb{P}^1$, \cite[Section
5.4]{lusnot}} \label{eqp1} Fix $k$ such that $1\leq k\leq n-1$ and
suppose $T = \CP^*\times \CP^*$ acts on $\mathbb{P}^1$ by
$(\lambda, \mu)\cdot (a: b) = ((\lambda\mu)^k a:
(\lambda^{-1}\mu)^{n-k} b)$. The following lemma can be proved by
explicit calculation on the two standard charts of $\mathbb{P}^1$.
\begin{lemma}
Let $T$ act on $\mathbb{P}^1$ as above. For every collection on
integers $i,i',j,j'$ satisfying $i'-j'=nm$ and $i-j=(2k-n)m$ for
some integer $m$, there exists a $T$-equivariant line bundle on
$\mathbb{P}^1$, unique up to isomorphism, where $(\lambda, \mu)$
acts as $\lambda^{j'}\mu^j$ on the fibre above $(1:0)$ and as
$\lambda^{i'}\mu^i$ on the fibre above $(0:1)$. All such
$T$-equivariant line bundles arise in this way.
\end{lemma}
We will denote the above line bundle $\mathcal{O}_k^{j'j,;i',i}$.

\subsection{Braid groups}
\label{brdgp}
The {\it affine Braid group} of type $\tilde{A}_{n-1}$, denoted $B_{ad}$, is the group with generators $\tilde{T}_i$ for $0\leq i\leq n-1$, satisfying the braid relations $\tilde{T}_i\tilde{T}_{i+1}\tilde{T}_i = \tilde{T}_{i+1}\tilde{T}_i\tilde{T}_{i+1}$ (we set $\tilde{T}_{n} = \tilde{T}_0$). There is a natural action of $<\sigma > = \mathbb{Z}/n\mathbb{Z}$ on $B_{ad}$ given by $\sigma (\tilde{T}_i) = \tilde{T}_{i+1}$. The \textit{extended affine Braid group} is
\[
B := \mathbb{Z}/n\mathbb{Z} \ltimes B_{ad}.
\]
Let $\iota : B \ra B$ be the involution which sends $\tilde{T}_i$
to $\tilde{T}_i^{-1}$ and fixes $\sigma$.

\subsection{Hecke algebras}
Let $\mathcal{A} = \mathbb{Z}[v^{\pm 1}]$. Let $\mathcal{A}[B]$ be the group algebra of $B$ over $\mathcal{A}$ and let $\mathcal{H}$,
the \textit{extended affine Hecke algebra}, be the quotient of $\mathcal{A}[B]$ by the ideal generated by $(\tilde{T}_i+v^{-1})(\tilde{T}_i - v)$
for $1\leq i\leq n$. The subalgebra generated by $T_i$ for $1\leq i\leq n-1$ (respectively $1\leq i\leq n$)
is the \textit{finite Hecke algebra} (respectively {\it affine Hecke algebra}) and denoted by $\mathcal{H}_{fin}$ (respectively $\mathcal{H}_{ad}$).

The involution $\iota$ on $B$ induces a ring automorphism on $\mathcal{H}$, which we denote by $\overline{\rule{0mm}{2mm}\;\;}$, sending $v$ to $v^{-1}$, $\tilde{T}_i$ to $\tilde{T}_i^{-1} = \tilde{T}_i + (v^{-1}-v)$ and fixes $\sigma$.

There is a second presentation of $\mathcal{H}$, \cite{lusaff}, discovered by Bernstein. In this presentation $\mathcal{H}$ is the $\mathcal{A}$-algebra generated by $\tilde{T}_i$, for $1\leq i\leq n-1$, and $\theta_x$, for $x\in P$, the weight lattice of $SL_n(\CP)$. For the fundamental weight $\varpi_{n-1}$ we have $\theta_{\varpi_{n-1}} = \sigma \tilde{T}_{n-2}\tilde{T}_{n-3}\ldots \tilde{T}_1\tilde{T}_n = \tilde{T}_{n-1}\tilde{T}_{n-2}\ldots \tilde{T}_2\tilde{T}_1\sigma$.

\subsection{}
\label{lusact} Since $T= \CP^*\times \CP^*$ we can identify
$\Rep(T)$ with $\mathbb{Z}[{v'}^{\pm 1}, v^{\pm 1}]$. The torus $T$
acts on $\mathcal{B}_{\chi}$ and $\Lambda_{\chi}$,
\cite[7.5]{slo}. The $T$-action on $\Pi_k$ agrees with the $k$-th
action considered in \ref{eqp1}, \cite[5.1]{lusnot}. There is an
action of $\mathcal{H}$ on $K_T(\mathcal{B}_{\chi})$,
\cite[Section 10]{lusbas1}. We will denote this action by
${\scriptscriptstyle \bullet}$ always. The action of $\mathcal{A}$
arises from the action of $\Rep(T)$, whilst the action of
$\tilde{T}_i$ for $1\leq i\leq n-1$ is explained in \cite[Sections
2,3]{lusnot}. For $x\in X$ the action of $\theta_x$ is given by
tensoring by certain line bundles on $\mathcal{B}$, \cite[Section
2]{lusnot}.
\subsection{}
\label{invol}
There is an involution,
\[
\tilde{\beta}: K_T(\mathcal{B}_{\chi}) \longrightarrow
K_T(\mathcal{B}_{\chi}),
\]
which is $\mathbb{Z}[{v'}^{\pm 1}]$-linear and twists the action of $\mathcal{H}$ by the involution $\overline{\rule{0mm}{2mm}\;\;}$, \cite[Proposition 12.10]{lusbas1}.

Let $\rule{0mm}{2mm}^{\dagger}$ denote the $\mathcal{A}$-algebra
involution of $\mathcal{A}[{v'}^{\pm 1}]$ which sends $v'$ to
${v'}^{-1}$ and let $\delta$ denote the $\mathcal{A}$-algebra
anti-involution of $\mathcal{H}$ which fixes $\tilde{T}_i$ and
sends $\sigma$ to $\sigma^{-1}$. By \cite[12.16]{lusbas1} there is
a pairing
\[
(\; | \; ): K_T(\mathcal{B}_{\chi})\times K_T(\mathcal{B}_{\chi})
\ra \mathbb{Z}[{v'}^{\pm 1}, v^{\pm 1}],
\]
which, for $F,F'\in K_T(\mathcal{B}_{\chi})$, $p\in
\mathbb{Z}[{v'}^{\pm 1},v^{\pm 1}]$ and $\tilde{T}\in
\mathcal{H}$, satisfies
\begin{enumerate}
\item $(p.F|F') = (F| p^{\dagger}.F') = p (F|F')$;
\item $(\tilde{T}{\scriptscriptstyle\bullet}F|F') = (F| \delta(\tilde{T}){\scriptscriptstyle\bullet} F')$;
\item $(F|F') = (F'|F)^{\dagger}$.
\end{enumerate}

As in \cite[5.11]{lusbas2} we define
\[
{\bf B}^{\pm}_{\mathcal{B}_{\chi}} = \{ F\in
K_T(\mathcal{B}_{\chi}) : \tilde{\beta}(F) = F, \; (F|F)\in 1+
\mathbb{Z}[{v'}^{\pm 1}, v^{-1}] \}.
\]
Thanks to \cite[Section 5]{lusnot} this set is a signed basis for
the free $\mathcal{A}$-module $K_T(\mathcal{B}_{\chi})$.

Let ${\bf p}_{k-1,k}$ denote the element in
$K_T(\mathcal{B}_{\chi})$ representing the skyscraper sheaf at
$p_{k-1,k}$ with trivial $T$-action. Abusing notation we will let
$\mathcal{O}_k^{j',j;i',i}$ denote the coherent sheaf on
$\mathcal{B}_{\chi}$ obtained by extension by zero of the
$T$-equivariant line bundle $\mathcal{O}_k^{j',j;i',i}$ on $\Pi_k
\cong \mathbb{P}^1$. For $1\leq k \leq n-1$ set ${\bf
O}_k=\mathcal{O}_k^{0,-n+k;-n,-k}$. We define
\[
{\bf O}_n = {\bf p}_{0,1} - \sum_{k=1}^{n-1}v^{n-k}{\bf O}_k = {\bf p}_{n-1,n} - \sum_{k=1}^{n-1}{v'}^nv^k {\bf O}_k.
\]
Then, by \cite[Proposition 5.25]{lusnot}, we have
\[
\textbf{B}^{\pm}_{\mathcal{B}_{\chi}} = \{ \pm {v'}^s {\bf O}_k :
1\leq k\leq n, \; s\in\mathbb{Z} \}.
\]
\subsection{}
\label{exp} There is an explicit description of the action of
$\mathcal{H}_{fin}$ on $K_T(\mathcal{B}_{\chi})$ given in
\cite[5.11]{lusnot}. For $1\leq i,k \leq n-1$ we have
\[
\tilde{T}_i {\scriptscriptstyle\bullet}{\bf O}_k = \begin{cases}
v{\bf O}_k \quad &\text{if $i = k$} \\
-v^{-1}{\bf O}_k - {\bf O}_{k\pm 1} &\text{if $i=k\pm 1$} \\
-v^{-1}{\bf O}_k &\text{otherwise,}
\end{cases}
\]
whilst $\tilde{T}_i({\bf p}_{0,1}) = -v^{-1}{\bf p}_{0,1} +
\delta_{1,i} (-{v'}^n + v^n){\bf O}_1$. By \cite[Lemma
5.24]{lusnot}, for $i\geq j$, the pairing is given by
\begin{equation}
\label{lprg} ({\bf O}_i,{\bf O}_j) = \begin{cases}
 1 + v^{-2} \quad &\text{if $i=j$} \\
-{v'}^{n}v^{-1} &\text{ if $i=n$ and $j=1$}\\
v^{-1} &\text{if $j=i-1$}\\
0 &\text{otherwise}.
\end{cases}
\end{equation}

\subsection{The equivariant Hilbert scheme}
Let $\Gamma$ be the finite cyclic subgroup of $SL_2(\CP)$ of order $n$, defined analogously to \ref{ks}. Let $\bf{H}$ be the $\Gamma$-equivariant Hilbert scheme of $\CP^2$, \cite[Chapter 4]{nak}. By definition the points of $\bf{H}$ are $\Gamma$-equivariant ideals of $\CP [X,Y]$, $I$, such that $\CP[X,Y]/I$ is isomorphic to the regular representation of $\Gamma$. There is a morphism
\[
\pi : {\bf H} \longrightarrow \CP^2/\Gamma,
\]
which, by \cite[Theorem 4.1]{nak}, is the minimal resolution of singularities of $\CP^2/\Gamma$. Let ${\bf H}_0$ denote the zero fibre $\pi^{-1}(0)$.
There is a $T$-equivariant isomorphism between ${\bf H}$ and $\Lambda_{\chi}$ which restricts to an isomorphism between ${\bf H}_0$ and $\mathcal{B}_{\chi}$.
\subsection{Derived equivalences}
Let $\mathcal{Z}={\bf H} \times_{\CP^2/\Gamma} \CP^2$. By
definition we have a commutative diagram
\[
\begin{CD}
\mathcal{Z} @> \psi >> & \CP ^2 \\
@V \eta VV & @VV \sigma V \\
\bf{H} @> \pi >> & \CP^2/\Gamma,
\end{CD}
\]
in which $\sigma$ and $\eta$ are finite, $\psi$ and $\pi$ are proper and birational, and $\eta$ is flat. Thanks to \cite[Theorem 1.4]{K-V} the functor
\[
\oo{\Phi} = ({\bf R}p_{\ast} \circ {\bf L}q^{\ast})^{\Gamma} : D_{\Gamma}(\CP ^2) \longrightarrow D(\bf{H})
\]
is an equivalence of triangulated categories. Moreover $\oo{\Phi}$ restricts to an equivalence between $D_{\Gamma}(\CP^2 , 0)$ and $D({\bf H}, {\bf H}_0)$, \cite[9.1]{bkr}.

Since $\psi$ and $\eta$ are $T$-equivariant morphisms $\oo{\Phi}$
can be lifted to an equivalence of triangulated categories
\[
\Phi = ({\bf R}p_{\ast} \circ {\bf L}q^{\ast})^{\Gamma} : D_{\Gamma\times T}(\CP^2) \longrightarrow D_T(\bf{H}).
\]
\subsection{}
\label{Tact}
Let $X$ be any $T$-equivariant quasi-projective variety and let $p_X: X\ra \{pt\}$ be projection onto a point. Let $v'$ (respectively $v$) be the one dimensional $T$-module where $(\lambda , \mu).1 = \lambda$ (respectively $(\lambda, \mu).1 = \mu$). As in \ref{easybun}, pulling back yields a $T$-equivariant line bundle which we denote by $v'_X$ (respectively $v_X$). Observe that tensoring by these line bundles realises an action of $\mathbb{Z}\times \mathbb{Z}$ on $\Coh_T(X)$ and $D_T(X)$, making them into (triangulated) $\mathbb{Z}\times \mathbb{Z}$-categories. We will denote the shift functors by ${v'}^i_{X}v_X^j$ for $i,j\in \mathbb{Z}$. When there is no confusion as to which variety we mean we will suppress the subscript.

In particular the above discussion applies to $\CP^2, {\bf H}, {\bf H}_0$ and $\mathcal{Z}$.
\begin{lemma}
The equivalence $\Phi$ is a $\mathbb{Z}\times \mathbb{Z}$-equivalence.
\end{lemma}
\begin{proof}
Let $\mathcal{F}\in D_{\Gamma\times T}(\CP^2)$. We have the following natural isomorphisms
\begin{eqnarray*}
{\bf R}p_{\ast} \circ ({\bf L}q^{\ast} ( - \otimes v_{\CP^2}))
&
\cong {\bf R}p_{\ast} \circ ({\bf L}q^{\ast}(-) \otimes q^{\ast}v_{\CP^2})
&
\cong {\bf R}p_{\ast} \circ ({\bf L}q^{\ast}(-) \otimes v_{\mathcal{Z}})
\\ &
\cong
{\bf R}p_{\ast} \circ ( {\bf L}q^{\ast} (- ) \otimes p^{\ast}v_{\bf{H}} )
&
\cong ({\bf R}p_{\ast}\circ {\bf L}q^{\ast} (-)) \otimes v_{\bf{H}},
\end{eqnarray*}
where the last isomorphism is the projection formula.
The same equations hold for $v'$, so $\Phi$ is a
$\mathbb{Z}\times\mathbb{Z}$-functor. Using the same formalism, one checks the inverse is also a $\mathbb{Z}\times\mathbb{Z}$-functor.
\end{proof}
\subsection{Serre-Grothendieck duality}
\label{sergro}
Let $\omega_{\CP^2}$ (respectively $\omega_{\bf{H}}$) be the canonical line bundle of $\CP^2$ (respectively $\bf{H}$).
As $\Gamma\times T$-equivariant bundles we have
$
\omega_{\CP^2} \cong v_{\CP^2}^{2}$.
As $T$-equivariant) bundles we have
$\omega_{\bf{H}} \cong v_{\bf{H}}^2$
\cite[Proposition 11.10]{lusbas1}.
We have contravariant isomorphisms of the derived categories $D_{\CP^2}$ and $D_{\bf{H}}$ where
\[
D_{\CP^2} = {\bf R} \mathcal{H}om_{\CP^2}( - , \omega_{\CP^2}[2]), \qquad D_{\bf{H}} = {\bf R} \mathcal{H}om_{\bf{H}}( - , \omega_{\bf{H}}[2]).
\]
Here $\mathcal{H}om$ denotes the sheaf of homomorphisms. Note that $D_{{\bf H}}$ sends $D_T({\bf H}, {\bf H}_0)$ to itself.
\subsection{Skew group rings}
\label{sgr}
The algebra $R=\CP [X,Y]$ has a $\mathbb{Z}\times \mathbb{Z}$-grading with $\text{deg}(X) = (1,1)$ and $\text{deg}(Y) = (-1,1)$. Since $\Gamma$ acts on $R$ we can form the skew group ring $R\ast \Gamma$. The $\mathbb{Z}\times \mathbb{Z}$-grading can be extended to $R\ast \Gamma$ by giving the elements of $\Gamma$ degree 0. We consider the category of finitely generated, bigraded $R\ast \Gamma$-modules, denoted $R\ast \Gamma\gmd$. As in \ref{grrg} $R\ast\Gamma\gmd$ is a $\mathbb{Z}\times \mathbb{Z}$-category. As in \ref{Tact}, for $i,j\in\mathbb{Z}$ we will denote the associated shift functor by ${v'}^iv^j$. The following lemma is standard.
\begin{lemma}
Taking global sections induces a $\mathbb{Z}\times \mathbb{Z}$-equivalence of categories
\[
\Upsilon: R\ast \Gamma\gmd \longrightarrow \Coh_{\Gamma \times T}(\CP^2).
\]
\end{lemma}
\subsection{}
\label{ye}
Let $\tau : R\ast \Gamma \rightarrow R\ast\Gamma$ be the anti-involution fixing $R$ and sending $g\in \Gamma$ to $g^{-1}$. If $M$ is a bigraded $R\ast \Gamma$-module then $\tau$ ensures $\Hom_R(M,R)$ is too. Serre-Grothendieck duality then induces the contravariant equivalence $D_{\CP^2} = {\bf R}\Hom_R( - , v^2R[2])$ on $D(R\ast\Gamma\gmd)$, the bounded derived category of $R\ast\Gamma\gmd$.

Let $D_0(R\ast\Gamma\gmd)$ be the bounded derived category of finitely generated, bigraded $R\ast\Gamma$-modules which are $(X,Y)$-primary, that is it which are annihilated by some power of $(X,Y)$. If $M$ is a complex in $D_0(R\ast\Gamma\gmd)$ then $D_{\CP^2}(M)\cong M^*$, the vector space dual of $M$, \cite[(4.9)]{yek}.
\subsection{}
\label{equi}
The results of this section can be summarised as follows.
\begin{prop}
There exists a $\mathbb{Z}\times \mathbb{Z}$-equivalence of triangulated categories
\[
\Psi = D_{\bf{H}} \circ \Phi \circ D_{\CP^2} \circ \Upsilon : D(R\ast \Gamma\gmd) \longrightarrow D_{T}(\bf{H}),
\]
which restricts to a $\mathbb{Z}\times \mathbb{Z}$-equivalence between $D_0(R\ast\Gamma\gmd)$ and $D_T({\bf H}_0)$.
\end{prop}
\section{The equivalence on $K$-theory}
\subsection{}
Let $K_0(R\ast\Gamma\gmd)$ be the Grothendieck group of $D_0(R\ast\Gamma\gmd)$. The equivalence of Proposition \ref{equi} yields isomorphisms
\[
\Sigma : K(R\ast\Gamma\gmd) \longrightarrow K_T({\bf H}) \quad \text{and}\quad \Sigma_0: K_0(R\ast\Gamma\gmd) \longrightarrow K_T({\bf H}_0).
\]
Moreover, since the equivalence preserves
$\mathbb{Z}\times \mathbb{Z}$-action,
both $\Sigma$ and $\Sigma_0$ are $\mathbb{Z}[{v'}^{\pm 1}, v^{\pm 1}]$-module isomorphisms.
\subsection{}
The simple $\Gamma$-modules are labelled by the elements of $\mathbb{Z}/n\mathbb{Z}$: the element $g$ acts on the $i$-th simple module as scalar multiplication by $\zeta^{i}$. Associated to the $i$-th simple $\Gamma$-module there are two bigraded $R\ast \Gamma$-modules:
the one dimensional module $S_i$ which is annihilated by both $X$ and $Y$, and the $R$-projective cover of $S_i$, denoted $R_i$.
As $\Gamma\times T$-equivariant sheaves on $\CP^2$ these correspond to the skyscraper sheaf supported at $0$
and the trivial line bundle, with trivial $T$-structure and $\Gamma$-structure given by the $i$-th simple $\Gamma$-module.
There is a Koszul resolution relating the two types of module
\begin{equation}
\label{Koszul}
0 \longrightarrow v^2R_i %\ar@{->}[r]^(.3){(X, -Y)^t}
\xrightarrow{(X, -Y)^t}
{v'}^{-1}vR_{i-1} \oplus v'vR_{i+1}
\xrightarrow{(Y, X)}
 R_i \longrightarrow
S_i \longrightarrow 0.
\end{equation}
By \ref{sergro} we have $D(R_i) = v^2R_{n-i}[2]$ and $D(S_i) = S_{n-i}$.
\subsection{Tautological bundles on \bf{H}}
\label{taut}
Since the projection $\eta: \mathcal{Z} \ra \bf{H}$ is a flat and finite
$\Gamma\times T$-equivariant morphism, the pushforward
$\eta_* \mathcal{O}_{\mathcal{Z}}$ is naturally a $\Gamma\times T$-equivariant
locally free sheaf, denoted $\mathcal{E}$. The fibre of $\mathcal{E}$ above
the ideal $I\in \bf{H}$ is the vector space $R/I$, the regular representation as a $\Gamma$-module.
Now $\Gamma$-equivariance allows us to decompose $\mathcal{E}$ into a direct sum of $T$-equivariant bundles
\[
\mathcal{E} = \bigoplus_{1\leq i\leq n} \mathcal{E}_i.
\]
By definition, the fibre of $\mathcal{E}_i$ above $I\in \bf{H}$ is the $i$-th isotypic component of $R/I$. These are described in \cite[5.27]{lusnot}.
\begin{lemma}
The isomorphism $\Psi$ sends $R_i$ to $\mathcal{E}_i^{\vee}$, the dual of $\mathcal{E}_i$.
\end{lemma}
\begin{proof}
Let $\pi_{\CP^2}:{\bf H}\times \CP^2 \ra \CP^2$ be the projection
map. Since there is a natural isomorphism ${\bf L}\phi^{\ast}
\cong \mathcal{O}_{\mathcal{Z}}\otimes^{{\bf L}}
\pi_{\CP^2}^{\ast}$, we have
\begin{eqnarray*}
\Psi (R_i) &\cong & D_{\bf{H}}(({\bf R}\eta_{\ast} {\bf
L}\phi^{\ast}(v_{\CP^2}^2R_{n-i}[2]))^{\Gamma}) \cong D_{\bf{H}}
(({\bf R}\eta_{\ast}( \mathcal{O}_{\mathcal{Z}}\otimes^{\bf L}
\pi_{\CP^2}^{\ast}( v_{\CP^2}^2R_{n-i}[2]))^{\Gamma})
\\
&\cong & D_{\bf{H}} ((v_{\bf H}^2\mathcal{E}[2]\otimes
S_{n-i})^{\Gamma}) \cong  D_{\bf{H}} (v_{\bf H}^2\mathcal{E}_i[2])
\cong \mathcal{E}_i^{\vee}.
\end{eqnarray*}
\end{proof}
\subsection{}
\label{mnth}
Using the resolution (\ref{Koszul}) and Lemma \ref{taut} we find that under $\Psi$ the module $S_i$ is sent to the Koszul complex
\begin{equation}
\label{kos}
\xymatrix{
v^2\mathcal{E}_i^{\vee} \ar@{->}[r] & {v'}^{-1}v\mathcal{E}_{i-1}^{\vee} \oplus v'v\mathcal{E}_{i+1}^{\vee} \ar@{->}[r] & \mathcal{E}_i^{\vee}.}
\end{equation}
\begin{prop}
For $1\leq i\leq n$ we have $\Sigma_0([S_i]) = {v'}^{n-i}{\bf O}_i$.
\end{prop}
\begin{proof}
In the non-equivariant setting the complex (\ref{kos}) was studied in \cite{itonak}. It is shown in \cite[Propositions 6.2]{itonak}
that its homology vanishes in degrees 1 and 2 and is $\mathcal{O}_{\Pi_i}(-1)$ in degree 0 if $1\leq i\leq n-1$.
Thus we have a quasi-isomorphism
\[
\xymatrix{
v^2\mathcal{E}_i^{\vee} \ar@{->}[r] \ar@{->}[d] & {v'}^{-1}v\mathcal{E}_{i-1}^{\vee} \oplus v'v\mathcal{E}_{i+1}^{\vee} \ar@{->}[r]^(.7){f} \ar@{->}[d] & \mathcal{E}_i^{\vee} \ar@{->}[d]\\
0 \ar@{->}[r] & 0 \ar@{->}[r] & \mathcal{E}_i^{\vee}/\im f.
}
\]
For $1\leq i\leq n-1$ the bundle $\mathcal{E}_i^{\vee}$ restricted to $\Pi_i$ is ${v'}^{n-i}{\bf O}_i$\cite[5.27]{lusnot}.
If $i=n$ a similar argument using \cite[6.4]{itonak}, shows that
$\Sigma_0 (S_i)= D_{\bf H}(\mathcal{O}_{{\bf H}_0})$.
It remains to prove that $[D_{\bf H}(\mathcal{O}_{{\bf H}_0})] = {\bf O}_n$.

For $1\leq k \leq n-1$ let $B_k$ be the subscheme of ${\bf H}_0$ consisting of the components $\Pi_1,\ldots ,\Pi_k$. Associated to this subscheme we have the maps of coherent sheaves arising from the inclusions $p_{k-1,k}\in \Pi_k \subseteq B_k$ and $p_{k-1,k}\in B_{k-1}\subset B_k$
\[
\mathcal{O}_{B_k} \longrightarrow \mathcal{O}_{B_{k-1}} \oplus \mathcal{O}_{\Pi_k} \longrightarrow p_{k-1,k}
\]
This is an exact sequence; so we have $[\mathcal{O}_{B_k}] = [\mathcal{O}_{B_{k-1}}] + [\mathcal{O}_{\Pi_k}] - {\bf p}_{k-1,k}$. Induction yields
\[
[\mathcal{O}_{{\bf H}_0}] = \sum_{k=1}^{n-1} [\mathcal{O}_{\Pi_k}] - \sum_{k=1}^{n-2} {\bf p}_{k,k+1}.
\]
By \cite[5.4]{lusnot} we have ${\bf p}_{k,k+1} = [\mathcal{O}_{\Pi_{k+1}}] - [\mathcal{O}_{k+1}^{n,-n+2(k+1);0,0}]$ and $[\mathcal{O}_{\Pi_1}] = {\bf p}_{0,1} + [\mathcal{O}_1^{n,-n+2;0,0}]$, which yields
\[
[\mathcal{O}_{{\bf H}_0}] = {\bf p}_{0,1} + \sum_{k=1}^{n-1}
[\mathcal{O}_k^{n,-n+2k;0,0}] = {\bf p}_{0,1} +
\sum_{k=1}^{n-1}{v'}^{n}v^{k} {\bf O}_k.
\]
By \cite[Lemma 5.16]{lusnot} we have $[D_{\bf{H}}({\bf{O}}_k) ] = -{v'}^nv^n
[{\bf{O}}_k]$ and $[D_{\bf{H}}(p_{0,1})] = p_{0,1}$. The result follows.
\end{proof}
The proposition shows that the signed basis $\{ \pm {v'}^s[S_i] :
1\leq i\leq n, s\in \mathbb{Z} \}$ of the free
$\mathcal{A}$-module $K_0(R\ast \Gamma\gmd)$ is sent by $\Sigma_0$
to ${\bf B}^{\pm}_{\mathcal{B}_\chi}$. A similar statement is true
for $\Lambda_{\chi}$, using $\Sigma$ and Lemma \ref{taut}.

\section{Categorification}
\label{catego}

\subsection{Braid group action}
\label{bga} Given a pair of objects $M,N\in D_0(R\ast\Gamma\gmd)$,
the functor ${\bf R}\Hom_{R\ast\Gamma}(M,N)$ takes values in the
derived category of bigraded, finite dimensional $\CP$-vector
spaces. We denote this category by $\mathcal{BV}$. In particular,
given the simple $R\ast\Gamma$-modules $S_i$ and $S_j$ the
resolution (\ref{Koszul}) shows
\begin{equation}
\label{ext}
{\bf R}^n \Hom_{R\ast\Gamma}(S_i,S_j) = \text{Ext}_{R\ast \Gamma}^n(S_i, S_j) =
\begin{cases}
\CP \quad &\text{if $n=0$ and $i=j$}\\
{v'}^{\pm 1}v^{-1}\CP  &\text{if $n=1$ and $i=j\pm 1$} \\
v^{-2}\CP  &\text{if $n=2$ and $i=j$}\\
0 &\text{otherwise}.
\end{cases}
\end{equation}
Given any complex $C$ belonging to $\mathcal{BV}$ and $M\in
D(R\ast\Gamma\gmd)$ we can form the tensor product $C\otimes D\in
D(R\ast\Gamma\gmd)$ and the space of linear maps $lin(C, M)\in
D(R\ast\Gamma\gmd)$, \cite[2a]{seitho}. Given a complex of
complexes $\cdots \ra C^{-1}\ra C^0 \ra C^1 \ra \cdots$ we write
$$ \{\cdots \ra C^{-1}\ra C^0 \ra C^1 \ra \cdots\} $$ to denote
the associated total complex.

For $0\leq i\leq n-1$ let $\tau_i : D_0(R\ast \Gamma\gmd)\ra
D_0(R\ast \Gamma\gmd)$ be the $\mathbb{Z}\times
\mathbb{Z}$-functor defined on objects by
%which associates to $M\in
%D_0(R\ast\Gamma\gmd)$ the complex
\[
\tau_i (M) =
\xymatrix{
\{{\bf R}\Hom (S_i, M) \otimes S_i \ar@{->}[r]^(0.7){ev} &M\},}
\]
where $ev$ is the evaluation map and $M$ is in degree $0$,
\cite[Definition 2.5]{seitho}.
By construction, this is nothing but the mapping cone of $ev$,
which turns out to be canonical for such complexes.
%Similarly, let $\tau_i^{-1}:D_0(R\ast \Gamma\gmd)\ra D_0(R\ast \Gamma\gmd)$ be
%the $\mathbb{Z}\times \mathbb{Z}$-functor which associates
%to $M\in D_0(R\ast\Gamma\gmd)$ the complex
%\[
%\xymatrix{
%\{M \ar@{->}[r]^(0.25){ev'} & lin({\bf R}\Hom_{R\ast \Gamma}(M, S_i),S_i) \},
%}\]
%where $ev'$ is the coevaluation map and $M$ is in degree $0$,
%\cite[Definition 2.7]{seitho}.
%Then $T_i$ and $T_i^{-1}$ are mutually inverse and
The various $\tau_i$'s furnish an action of the affine braid group,
$B_{ad}$, on $D_0(R\ast \Gamma\gmd)$, \cite[Proposition 2.4, 2c and Example 3.9]{seitho}. For our purposes it is better to work
with an adjusted action for $0\leq i\leq n-1$
\[
\tilde{\tau}_i = D_{\CP^2} \circ \tau_{n-i} \circ D_{\CP^2} \circ
v^{-1} \circ [1].
\]
The $\mathbb{Z}\times \mathbb{Z}$-equivalences $\tilde{\tau}_i$
continue to satisfy the braid relations.
\subsection{}
\label{brdcalc} For $i\in \mathbb{Z}/n\mathbb{Z}$ let $M^+_i$
(respectively $M_i^-$) be the unique two dimensional bigraded
$R\ast\Gamma$-module with head isomorphic to $S_i$ and socle
isomorphic to $v'vS_{i+1}$ (respectively ${v'}^{-1}vS_{i-1}$),
where multiplication by $X$ (respectively $Y$) sends the head to
the socle.
\begin{lemma}
For $i,j\in\mathbb{Z}/n\mathbb{Z}$ we have
\[
\tilde{\tau}_i(S_j) =
\begin{cases}
vS_j \quad &\text{if $i=j$} \\
v^{-1}M_j^{\pm}[1] &\text{if $i = j\pm 1$} \\
v^{-1}j[1] &\text{otherwise}.
\end{cases}
\]
\end{lemma}
\begin{proof}
We first calculate $\tau_i$. If $i=j$ then, using (\ref{Koszul}),
we see that $\tau_j(S_j)$ equals the complex $$ v^2R_j
\xrightarrow{(X, -Y)^t} {v'}^{-1}vR_{j-1} \oplus {v'}vR_{j+1}
\xrightarrow{{\left( \begin{matrix} \scriptstyle{Y} &
\scriptstyle{0}
\\ \scriptstyle{X} & \scriptstyle{0} \\
\end{matrix}
\right)}}
 R_j \oplus R_j
\xrightarrow{\left(\begin{matrix} \scriptstyle{\epsilon_j} &
\scriptstyle{0} & \scriptstyle{0} \\ \scriptstyle{\epsilon_j} &
\scriptstyle{X} & \scriptstyle{-Y}
\end{matrix} \right)^t}
$$ $$ \ra  S_i\oplus {v'}^{-1}v^{-1}R_{j-1}\oplus
{v'}v^{-1}R_{j+1} \xrightarrow{(0, Y, X)} v^{-2}R_j. $$ It is
straightforward to check that this complex only has homology at
its end term, equal to $v^{-2}S_j$. Therefore $\tau_j(S_j) =
v^{-2}S_j[1]$. If $i=j+1$ (the case $i=j-1$ is analogous) we find
$\tau_{j+1}(S_j)$ is the complex $$ {v'}v R_{j+1}
\xrightarrow{(-X, Y)^t} R_{j} \oplus {v'}^2R_{j+2}
\xrightarrow{\left( \begin{matrix} \epsilon_j & 0 \\ -Y & -X
\end{matrix} \right)} S_j \oplus {v'}v^{-1}R_{j+1}. $$ Again, it is
straightforward to check that this complex only has homology at
its end term. A basis for this homology group is $\{ \oo{(1,0)},
\oo{(0,1)}\}\subset S_j \oplus v'v^{-1}R_{j+1}$, showing that this
is the unique bigraded $R\ast\Gamma$-module with head
$v'v^{-1}S_{j+1}$ and socle $S_j$ (multiplication by $Y$ sends the
head $\oo{(0,1)}$ to the socle $\oo{(1,0)}$). Call this $N_j$. If
$i\neq j,j\pm 1$, then $\tau_i(S_j) = S_j$ since the relevant
Ext-group vanishes, (\ref{ext}).

The description of $\tilde{\tau}_i$ follows immediately, noting
that $D_{\CP^2}(N_{n-j}) \cong M^-_j$ thanks to \ref{ye}.
\end{proof}
\subsection{}
\label{cycle} Let $\varsigma$ be the automorphism of $R\ast
\Gamma$ which fixes $X$ and $Y$ and sends $g$ to $\zeta g$.
Following \ref{grrg}, this induces an $\mathbb{Z}\times
\mathbb{Z}$-equivariant action of $\mathbb{Z}/n\mathbb{Z}$ on the
derived category $D(R\ast \Gamma\gmd)$. This satisfies
$\varsigma\circ \tilde{\tau}_i \circ \varsigma^{-1} \cong
\tilde{\tau}_{i+1}$. Therefore we have an action of the extended
affine braid group $B$ on $D(R\ast \Gamma\gmd)$ via
$[\tilde{T}_i]=\tilde{\tau}_i$, $[\sigma ]= \varsigma$.

\subsection{}
\label{agree}
The action of $B$ on $D_0(R\ast\Gamma\gmd)$ induces an action of $\mathcal{A}[B]$ on $K_0(R\ast\Gamma\gmd)$ which, by Lemma \ref{brdcalc}, satisfies
\[
\tilde{T}_i([S_j]) = \begin{cases}
 v[S_j] \quad & \text{if $i=j$}\\
-v^{-1}[S_j] - {v'}^{\pm 1}[S_{j\pm 1}] &\text{if $i=j\pm 1$} \\
-v^{-1}[S_j] &\text{otherwise}.
\end{cases}
\]
In particular, since $(\tilde{T}_i +v^{-1})(\tilde{T}_i - v)$ acts as zero, this factors through an action of $\mathcal{H}$.
\begin{prop}
The isomorphism $\Sigma_0 : K_0(R\ast\Gamma\gmd) \ra K_T({\bf H}_0)$ is $\mathcal{H}_{ad}$-equivariant.
\end{prop}
The action on $K_T({\bf H}_0)$ referred to in the proposition
comes from \ref{lusact}. The proposition will be proved in the
following six subsections, where we also identify the full
$\mathcal{H}$-action.
\subsection{}
Let's translate the action from $K_0(R\ast\Gamma\gmd)$ over to $K_T({\bf H}_0)$. Using Proposition \ref{mnth} it is straightforward to check we have
\[
\tilde{T}_i({\bf O}_j) =
\begin{cases}
v{\bf O}_j \quad &\text{if $i=j$} \\
-v^{-1}{\bf O}_1 - {v'}^{-n}{\bf O}_n &\text{if $i=n, j=1$} \\
-v^{-1}{\bf O}_n - {v'}^n{\bf O}_1 &\text{if $i=1, j=n$} \\
-v^{-1}{\bf O}_j - {\bf O}_{j\pm 1} &\text{if $i=j\pm 1$ is not as above} \\
-v^{-1}{\bf O}_j &\text{otherwise}.
\end{cases}
\]
By \ref{exp} the action of $\mathcal{H}_{fin}$ here agrees with the ${\scriptscriptstyle \bullet}$-action on ${\bf O}_j$ for $1\leq j\leq n-1$.
\subsection{}
We first check that the actions of $\mathcal{H}_{fin}$ also agree
on ${\bf O}_n$. Thanks to \ref{exp} we have $$
\tilde{T}_1{\scriptscriptstyle \bullet}{\bf O}_n = -v^{-1}{\bf
p}_{0,1} + (-{v'}^n + v^n){\bf O}_1 - v^n{\bf O}_1 -
v^{n-2}(-v^{-1}{\bf O}_2 - {\bf O}_1) + \sum_{k=3}^{n-1}
v^{n-k-1}{\bf O}_k $$ $$ =-v^{-1}{\bf p}_{0,1} + (-{v'}^n +
v^{n-2}){\bf O}_1 + \sum_{k=2}^{n-1} v^{n-k-1}{\bf O}_k =
-v^{-1}{\bf O}_n - {v'}^n{\bf O}_1. $$ Similarly, for $2\leq i\leq
n-2$, we have $\tilde{T}_i{\scriptscriptstyle \bullet}{\bf O}_n
=
% -v^{-1}{\bf p}_{0,1} + \sum_{k=1}^{n-1} v^{n-k-1}{\bf O}_k + v^{n-i+1}(v^{-1}{\bf O}_{i-1}+{\bf O}_i) - v^{n-i+1}{\bf O}_i + v^{n-i-1}(v^{-1}{\bf O}_{i+1} + {\bf O}_i)
-v^{-1}{\bf O}_n$, whilst $\tilde{T}_{n-1}{\scriptscriptstyle
\bullet}{\bf O}_n =
%-v^{-1}{\bf p}_{0,1} + \sum_{k=1}^{n-3}v^{n-k-1}{\bf O}_k + v^2(v^{-1}{\bf O}_{n-2} - {\bf O}_{n-1}) -v^2{\bf O}_{n-1}
-v^{-1}{\bf O}_n - {\bf O}_{n-1}$. This shows that $\Sigma_0$ is
$\mathcal{H}_{fin}$-equivariant.
\subsection{}
\label{lustwi} We now calculate the action of
$\theta_{\varpi_{n-1}}= \sigma
\tilde{T}_{n-2}\tilde{T}_{n-3}\ldots \tilde{T}_1\tilde{T}_n$ in
both cases. Tensoring by the line bundle $L_{n-1,n}$ on
$\mathcal{B}$, described in \cite[5.1]{lusnot}, corresponds to
$\theta_{\varpi_{n-1}}{\scriptscriptstyle \bullet}$. The
$T$-equivariant structure of $L_{n-1,n}$ is given in
\cite[5.2]{lusnot} for $1\leq k\leq n$ as follows
\[
{L}_{n-1,n}|_{p_{k-1,k}} =
\begin{cases}
{v'}^{n-1} \quad & \text{if $k=n$} \\
{v'}^{-1}v^{2-n} & \text{if $k < n$}.
\end{cases}
\]
In particular, if $k < n$ then the restriction of ${L}_{n-1,n}$ to
$\Pi_k$ is the trivial line bundle with equivariant structure
${v'}^{-1}v^{2-n}$. We deduce for $1\leq k \leq n-2$
\[
\theta_{\varpi_{n-1}}{\scriptscriptstyle \bullet} {\bf O}_k  =
{L}_{n-1,n}\otimes {\bf O}_k = {v'}^{-1}v^{2-n} {\bf O}_k.
\]
We have
\begin{eqnarray*}
\theta_{\varpi_{n-1}}{\scriptscriptstyle \bullet} {\bf O}_{n-1} &
= & {L}_{n-1,n} \otimes {\bf O}_{n-1} = [O_{n-1}^{-1,2-n; n-1,0}]
\otimes [O_{n-1}^{0,-1;-n,1-n}] \\ & = & [O_{n-1}^{-1,1-n;-1,1-n}]
=  {v'}^{-1}v^{1-n} [O_{n-1}^{0,0;0,0}] \\ &=& {v'}^{-1}v^{1-n}(
{\bf p}_{n-1,n} + [O_{n-1}^{0,0;-n,2-n}]) = {v'}^{-1}v^{1-n}({\bf
p}_{n-1,n} + v{\bf O}_{n-1}) \\ & = & {v'}^{-1}v^{1-n}({\bf
p}_{n-1,n} - \sum_{k=1}^{n-1}{v'}^nv^k{\bf O}_k) +
\sum_{k=1}^{n-2} {v'}^{n-1}v^{k+1-n}{\bf O}_k +
({v'}^{n-1}+{v'}^{-1}v^{2-n}){\bf O}_{n-1} \\ & =&
{v'}^{n-1}v^{1-n}(\sum_{k=1}^{n-1} v^i{\bf O}_i) +
{v'}^{-1}v^{2-n}{\bf O}_{n-1} + {v'}^{-1}v^{1-n}{\bf O}_n,
\end{eqnarray*}
and similarly
\begin{eqnarray*}
\theta_{\varpi_{n-1}}{\scriptscriptstyle \bullet} {\bf O}_{n} & =
& {L}_{n-1,n} \otimes ({\bf p}_{0,1} - \sum_{k=1}^{n-1}v^{n-k}{\bf
O}_k) = {v'}^{-1}v^{2-n}{\bf p}_{0,1} - \sum_{k=1}^{n-2}
{v'}^{-1}v^{2-k}{\bf O}_k - \\ & - &
v({v'}^{n-1}v^{1-n}(\sum_{k=1}^{n-1} v^i{\bf O}_i) +
{v'}^{-1}v^{2-n}{\bf O}_{n-1} + {v'}^{-1}v^{1-n}{\bf O}_n) \\ &=
&- \sum_{k=1}^{n-2}{v'}^{n-1}v^{2+k-n}{\bf O}_k - {v'}^{n-1}v{\bf
O}_{n-1} = -{v'}^{n-1}v^{2-n}(\sum_{k=1}^{n-2} v^k{\bf O}_k).
\end{eqnarray*}

\subsection{}
We calculate $\tilde{T}_{n-2}\tilde{T}_{n-3}\ldots \tilde{T}_1\tilde{T}_0({\bf O}_i)$ for $1\leq i\leq n$. Note that $\tilde{T}_i\tilde{T}_{i-1}({\bf O}_i) = v^{-1}{\bf O}_i$ for $2\leq i\leq n-1$, whilst $\tilde{T}_1\tilde{T}_n({\bf O}_1) = {v'}^{-n}v^{-1}{\bf O}_n$.
It follows that, for $2\leq i\leq n-1$,
\begin{equation}
\label{e1}
\tilde{T}_{n-2}\tilde{T}_{n-3}\ldots \tilde{T}_1\tilde{T}_0({\bf O}_i)= (-1)^{n-1}v^{2-n} {\bf O}_{i-1},
\end{equation}
whilst
\begin{equation}
\label{e2}
\tilde{T}_{n-2}\tilde{T}_{n-3}\ldots \tilde{T}_1\tilde{T}_0({\bf O}_1)= (-1)^{n-1}{v'}^{-n}v^{2-n}{\bf O}_n.
\end{equation}
We have $\tilde{T}_j \ldots \tilde{T}_1\tilde{T}_0 ({\bf O}_{n-1}) = (-1)^{j+1}(v^{-j-1}{\bf O}_{n-1} + v^{-j}{\bf O}_n + {v'}^nv^{-j}\sum_{k=1}^j v^k{\bf O}_k)$ for $1\leq j\leq n-3$: this can be proved by induction. We deduce that
$$
\tilde{T}_{n-2}\tilde{T}_{n-3}\ldots \tilde{T}_1\tilde{T}_0({\bf O}_{n-1})
\; = \;
\tilde{T}_{n-2}((-1)^{n-2}(v^{2-n}{\bf O}_{n-1} + v^{3-n}{\bf O}_n + {v'}^nv^{3-n}\sum_{k=1}^{n-3} v^k {\bf O}_k))
$$
\begin{equation}
\label{e3}
= (-1)^{n-1}(v^{1-n}{\bf O}_{n-1} + v^{2-n}{\bf O}_n + {v'}^nv^{2-n}\sum_{k=1}^{n-2} v^k{\bf O}_k + v^{2-n}{\bf O}_{n-2}).
\end{equation}
Arguing by induction, we see $\tilde{T}_j\ldots
\tilde{T}_1\tilde{T}_0 ({\bf O}_n) = (-1)^j(v^{1-j}{\bf O}_n +
{v'}^nv^{-j+1}\sum_{k=1}^j v^k {\bf O}_k)$ for $1\leq j\leq n-2$.
In particular
\begin{equation}
\label{e4}
\tilde{T}_{n-2}\tilde{T}_{n-3}\ldots \tilde{T}_1\tilde{T}_0({\bf O}_n)
= (-1)^{n-2}(v^{3-n}{\bf O}_n + {v'}^n v^{3-n}\sum_{k=1}^{n-2}v^k {\bf O}_k).
\end{equation}
\subsection{}
\label{kaptwi}
To get $\theta_{\varpi_{n-1}}$ we must twist by the automorphism $\sigma$. This twist sends $S_i$ to $S_{i+1}$ so, by Theorem \ref{mnth}, sends ${\bf O}_i$ to ${v'}^{-1}{\bf O}_{i+1}$ for $1\leq i\leq n-1$
and ${\bf O}_n$ to ${v'}^{n-1}{\bf O}_1$. Combining this with (\ref{e1}), (\ref{e2}), (\ref{e3}) and (\ref{e4}) yields
\[
\theta_{\varpi_{n-1}}({\bf O}_i) =
\begin{cases}
(-1)^{n-1}{v'}^{-1}v^{2-n} {\bf O}_i \quad &\text{if $1\leq i \leq n-2$} \\
(-1)^{n-1}({v'}^{n-1}v^{1-n}(\sum_{k=1}^{n-1} v^i{\bf O}_i) + {v'}^{-1}v^{2-n}{\bf O}_{n-1} + {v'}^{-1}v^{1-n}{\bf O}_n) &\text{if $i=n-1$}\\
(-1)^{n-2}{v'}^{n-1}v^{2-n}(\sum_{k=1}^{n-2} v^k{\bf O}_k) &\text{if $i=n$}.
\end{cases}
\]
\subsection{}
Comparing \ref{lustwi} and \ref{kaptwi} we see the actions of $\theta_{\varpi_{n-1}}$ differ by scalar multiplication by $(-1)^{n-1}$. Since $\theta_{\varpi_{n-1}} =\tilde{T}_{n-1}\tilde{T}_{n-2}\ldots \tilde{T}_1 \sigma$ and the actions of $\mathcal{H}_{fin}$ agree we deduce that the actions of $\sigma$ must differ by scalar multiplication by $(-1)^{n-1}$.
Since $T_0 = \sigma T_{n-1} \sigma^{-1}$ we deduce that the actions agree on $T_n$, and hence on $\mathcal{H}_{ad}$, proving Theorem \ref{agree}.

There is an involution of $\mathcal{H}$ which fixes $\mathcal{H}_{ad}$ and sends $\sigma$ to $(-1)^{n-1}\sigma$. The calculations show that the two module structures on $K_T({\bf H}_0)$ are twists of each other under this involution.
\subsection{Duality}
\label{dua}
Let $\nu$ be the automorphism of $R\ast \Gamma$ which swaps $X$ and $Y$ and sends $g$ to $g^{-1}$. As in \ref{cycle} we can twist an object $M\in R\ast \Gamma\gmd$ by $\nu$ where we set
\[
(M^{\nu})_{i,j} = M_{-i,j}
\]
for $i,j\in\mathbb{Z}$.
Twisting by $\nu$ commutes with $D_{\CP^2}$ by \ref{ye}. There is a contravariant self-equivalence
\[
\tilde{\beta} : D_0(R\ast \Gamma\gmd) \longrightarrow D_0(R\ast \Gamma\gmd)
\]
which sends $M$ to $D_{\CP^2}(M)^{\nu}$. Note that
$\tilde{\beta}({v'}^av^bM)= {v'}^av^{-b}\tilde{\beta}(M)$ for
$a,b\in\mathbb{Z}$ and $M\in D_0(R\ast\Gamma\gmd)$. Moreover,
$\tilde{\beta}(S_i) = S_i$ for $1\leq i\leq n$.

Recall the involution $\iota$ on $B$, defined in \ref{brdgp}.
\begin{lemma}
For all $b\in B$ and $M\in D_0(R\ast\Gamma\gmd)$ we have
$\tilde{\beta}(b(M))\cong\iota(b)(\tilde{\beta}(M))$.
\end{lemma}
\begin{proof}
Given any object $C\in \mathcal{BV}$ and $M\in D(R\ast\Gamma\gmd)$
evaluation gives a natural isomorphism $lin(C, D_{\CP^2}(M))
\longrightarrow D_{\CP^2}(C\otimes M)$ in $D(R\ast \Gamma\gmd)$.
Hence, by \cite[Section 2]{seitho},
\begin{eqnarray*}
T_i^{-1}(M) & = & \{ \xymatrix{M \ar@{->}[r]^(0.15){ev'} &lin ({\bf R}\Hom_{R\ast\Gamma} (M, D_{\CP^2}S_{n-i}), D_{\CP^2}S_{n-i})} \} \\
& \cong & \{ M \longrightarrow D_{\CP^2}({\bf R}\Hom_{R\ast\Gamma} (M, D_{\CP^2}S_{n-i})\otimes S_{n-i}) \}
\\
& \cong &  \{ M \longrightarrow D_{\CP^2}({\bf R}\Hom_{R\ast\Gamma} (S_{n-i}, D_{\CP^2}(M))\otimes S_{n-i}) \} \\
& \cong & D_{\CP^2}( \{ \xymatrix{({\bf R}\Hom_{R\ast\Gamma} (S_{n-i}, D_{\CP^2}(M))\otimes S_{n-i}) \ar@{->}[r]^(0.75){ev} & D_{\CP^2}(M) }\}) \\
& = & D_{\CP^2} T_{n-i}(D_{\CP^2}(M)).
\end{eqnarray*}
By construction $\rule{0mm}{2mm}^{\nu}\circ T_i = T_{n-i}\circ\rule{0mm}{2mm}^{\nu}$ so we deduce that $\tilde{\beta} \circ T_i = T_i^{-1} \circ \tilde{\beta}$. The same statement with $\tilde{T}_i$ follows since
\[
\tilde{T}_i\circ\tilde{\beta} = D_{\CP^2}\circ T_{n-i} \circ D_{\CP^2} \circ v^{-1}\circ [1]\circ \tilde{\beta}= \tilde{\beta}\circ D_{\CP^2} \circ T_{n-i}^{-1} \circ D_{\CP^2} \circ v\circ [-1] = \tilde{\beta}\circ \tilde{T}_i^{-1}.
\]
Finally, since $\nu\sigma = \sigma^{-1}\nu$ and $(M^*)^{\sigma} =
(M^{\sigma^{-1}})^*$ we find $(D_{\CP^2}(M^{\sigma}))^{\nu} =
(D_{\CP^2}(M))^{\sigma^{-1}\nu} = (D_{\CP^2}(M)^{\nu})^{\sigma}$,
as required.
\end{proof}

\subsection{A pairing}
\label{prg}
There is a pairing
\[
(\: , \,) : D(R\ast\Gamma \gmd )\times D(R\ast\Gamma \gmd) \longrightarrow \mathcal{BV},
\]
given by $(M, N) = {\bf R}\Hom_{R\ast \Gamma}(\tilde{\beta}(N), M)$.

Recall from \ref{invol} the involution $\rule{0mm}{2mm}^{\dagger}$ on $\mathbb{Z}[{v'}^{\pm 1}, v^{\pm 1}]$. This induces an involution on $\mathcal{BV}$ which we denote also by $\rule{0mm}{2mm}^{\dagger}$.
\begin{lemma}
The pairing satisfies the following properties:

1) $({v'}^av^bM,N) = (M, {v'}^{-a}v^bN) = {v'}^av^b(M,N)$;

2) $(T_i M , N) = (M , T_i N)$ and $(\sigma M, N) = (M, \sigma^{-1}N)$;

3) $(M, N) = (N, M)^{\dagger}$.
\end{lemma}
\begin{proof}
Part 1 follows immediately from the variance properties of ${\bf R}\Hom$ and $\tilde{\beta}$. By Lemma \ref{dua} we have
\[
{\bf R}\Hom_{R\ast\Gamma}(\tilde{\beta}(N), T_i(M) ) = {\bf R}\Hom_{R\ast\Gamma}( T_i^{-1}(\tilde{\beta}(N)),M) = {\bf R}\Hom_{R\ast\Gamma} ( \tilde{\beta}(T_i(N)),M),
\]
and
\[
{\bf R}\Hom_{R\ast\Gamma}(\tilde{\beta}(N),M^{\sigma}) = {\bf R}\Hom_R{\ast \Gamma}(\tilde{\beta}(N)^{\sigma^{-1}},M) = {\bf R}\Hom_{R\ast\Gamma} (\tilde{\beta}(N^{\sigma}),M),
\]
proving part 2. Finally, we have
\[
{\bf R}\Hom_{R\ast\Gamma}(\tilde{\beta}(N), M) \cong {\bf R}\Hom_{R\ast\Gamma}(D_{\CP^2}(M), D_{\CP^2}(\tilde{\beta}(N))) \cong {\bf R}\Hom_{R\ast\Gamma}(D_{\CP^2}(M), N^{\nu}).
\]
In $R\ast\Gamma\gmd$ the identity provides an identification
$\Hom_{R\ast\Gamma}(A,B) =
\Hom_{R\ast\Gamma}(A^{\nu},B^{\nu})^{\dagger}$, which, together
with above equation, proves part 3.
\end{proof}
\subsection{}
\label{cat} The next corollary shows that the structures defined
in this chapter give a categorification of the action of
$\mathcal{H}_{ad}$ on $K_T(\mathcal{B}_{\chi})$ discussed in
\ref{lusact}. It is clear we can extend this to a categorification
of the $\mathcal{H}$-action. For more details about
categorification, see \cite{bfk}.

\begin{theorem}
The duality $\tilde{\beta}$ and pairing $(\;,\;)$ on $D_0(R\ast
\Gamma\gmd)$ induce a duality and a pairing on
$K_0(R\ast\Gamma\gmd)$, which, under the $\mathcal{H}_{ad}$-module
isomorphism $\Sigma_0$, correspond to the duality and pairing
defined in \ref{invol}.
\end{theorem}
\begin{proof}
By construction $\tilde{\beta}({v'}^av^b[i^!])= {v'}^av^{-b}[i^!]$ for $a,b\in\mathbb{Z}$ and $1\leq i\leq n$. The claim about duality then follows from the sentence following the proof of Proposition \ref{mnth}.

The calculation of Ext-groups in (\ref{ext}) shows that
 \[
(S_i, S_j) = \begin{cases}
1 + v^{-2}[-2] \quad &\text{if $j=i$} \\
{v'}^{\pm 1}v^{-1}[-1] &\text{if $j= i\pm 1$} \\
0 &\text{otherwise}
\end{cases}
\]
Calculating the image of this in $K_T({\bf H}_0)$ under $\Sigma_0$ yields (\ref{lprg}). Thus the pairings agree since Claims 1,2 and 3 in Lemma \ref{prg} correspond to \ref{invol}.1,\ref{invol}.2 and \ref{invol}.3
\end{proof}

\section{Final Remarks}
\label{finalr}
\subsection{Coinvariant algebra}
Define the \textit{skew coinvariant algebra}
\[
\cl = \frac{\LF[X,Y]}{(X^n, XY, Y^n)}\ast \Gamma.
\]
Since $(X^n, XY, Y^n)$ is a homogeneous ideal for the bigrading on
$\LF[X,Y]$ similar to that in \ref{sgr}, there is a natural
bigraded structure on $\cl$. Let $\cl\gmd$ denote full
$\mathbb{Z}\times \mathbb{Z}$-subcategory of
$\LF[X,Y]\ast\Gamma\gmd$ consisting of $\cl$-modules.

Observe that the constructions in Section \ref{catego} of the
braid group action, duality and pairing carry over to $\LF
[X,Y]\ast\Gamma$ if the characteristic of $\LF$ is either zero or
coprime to $n$. We will use $\tilde{\beta}_{\LF}$ to denote the
duality. In particular, passing to $K$-theory induces an action of
$\mathcal{H}$ on the Grothendieck group of bigraded
$C_{\LF}(n)$-modules.

There is an isomorphism
\[
\upsilon: N_{\LF}(n) \longrightarrow \cl,
\]
which sends $X$ to $\sum_{k=1}^n b_k$, $Y$ to $\sum_{k=1}^n a_k$ and $g$ to $\sum_{k=1}^n \zeta^{n-k}e_k$. By construction the $\mathbb{Z}$-grading on
$N_\LF(n)$ corresponds to the $v'$-component on $\cl$.
We write $\cl\gmd|_{v=1}$ to denote to the associated category of
$\mathbb{Z}$-graded $\cl$-modules.

\subsection{}
\label{dua2} Fix a regular weight $\lambda\in C_0$. Proposition
\ref{lienoc} proves the existence of a $\mathbb{Z}$-equivalence of
categories
\[
\Theta : \uclt\md \longrightarrow C_\K (n) \gmd|_{v=1}.
\]
Under the isomorphism $\upsilon$, the element $g$ acts on $e_i$ by scalar multiplication by $\zeta^{n-i}$. Therefore $\Theta(L_i) = S_{n-i}$.
The reader can check that
the duality $D$ on $\uclt\md$
corresponds, through $\Theta$,
to $\tilde{\beta}_\K$.

\subsection{Translation functors (II)}
\label{tran2}
Let
$
\gamma_i = T_{\mu_i}^{\lambda}\circ T_{\lambda}^{\mu_i} : \BT\md
\longrightarrow \BT\md
$
be the $i$th \textit{wall-crossing functor}. By the adjunction
property, there is a natural transformation $\epsilon_i : \gamma_i
\ra \textit{id}$, which, as in \cite[Section 3]{ric}, induces a
triangulated $\mathbb{Z}$-functor $\tilde{T}_i$ on the bounded
derived category of $\BT\md$.

The functor $\tilde{T}_i$ does not necessarily  restrict to
$D(\uclt\md)$. However, there is an isomorphism between
$K(\BT\md)$ and $K(\uclt\md)$. By transfer of structure, the
$\mathbb{Z}[{v'}^{\pm 1}]$-linear operators $\tilde{T}_i$ act on
$K(\BT\md)$.

Recall the action of $\mathcal{H}$ on $K(C\gmd)$. This gives us
operators $\tilde{T}_i$ on $K(C\gmd|_{v=1})$ when we specialise
$v$ to $1$. It is not difficult to show that
\begin{center}
{\em The isomorphism $\theta: K(\uclt\md) \ra K(C\gmd|_{v=1})$
commutes with}
\\
{\em the operators $\tilde{T}_i$ for $0\leq i\leq n-1$}.
\end{center}
\subsection{Extension to type B}
For subregular representations of Lie algebras of type $B$,
Jantzen has proved an analogue of Proposition \ref{uni} and
calculated several extension groups, \cite[Section 3]{jan}. Then,
{\em if we can prove an analogue of Proposition \ref{pro}}, it
follows formally from the arguments of Sections \ref{red},
\ref{hq} and \ref{prec} that the central reduction of a block of a
subregular reduced enveloping algebra of type $B$ is a matrix ring
over a central reduction of Hodges' deformation of a Kleinian
singularity of type $A$.

Unfortunately, the proof of Proposition \ref{pro} does not
immediately generalise to type $B$ since it is no longer true that
we can find all highest weights in the fundamental alcove, $C_0$.
For regular weights, however, this can be remedied as follows. Let
$\bc$ denote a block of a subregular reduced enveloping algebra of
type $B_n$ associated to a regular weight $\lambda$. Rickard's
results in \cite[Section 3]{ric} remain valid in this situation.
To check this requires the specific information on the
wall-crossing functors, denoted $\gamma_{i}$ for $1\leq i\leq n$,
provided by Jantzen in \cite[Section H]{ja1}. Thus, for $1\leq
i\leq n$, there are derived self-equivalences $F_i$ on the bounded
derived categories of $\bc$-modules. In particular, given any
$\bc$-module $M$, $F_i(M)$ is the complex $\gamma_{i}(M)\ra M$,
the map being given by the counit of $\gamma_{i}$. Combining the
known behaviour of baby Verma modules under wall-crossing,
\cite[11.20]{ja3}, with existing results on filtrations of the
projective indecomposable modules by baby Verma modules,
\cite[Proposition 10.11]{ja3}, it can be shown that, given two
projective indecomposable $\bc$-modules $Q$ and $Q'$, there are
integers $1\leq i_1,\ldots ,i_t\leq n$ such that $F_{i_1}\ldots
F_{i_t}(Q)$ is quasi-isomorphic to $Q'$. Since it is known that
the centre of $\bc$ is a derived-invariant, the proof of
Proposition \ref{pro} can be generalised, using the derived
category, to deal with $\bc$. As a result we find the central
reduction of $\bc$ is a matrix ring over the central reduction of
Hodges' deformation of a Kleinian singularity of type $A_{2n-1}$.

\end{document}